\documentclass[leqno,11pt]{article}

\usepackage{latexsym}

\makeatletter
\@addtoreset{equation}{section}
\makeatother

\mathchardef\emptyset="001F

\newcommand{\R}{I\hspace{-1.5mm}R}
\newcommand{\scriptR}{I\hspace{-0.9mm}R}

\newtheorem{Theorem}{Theorem}[section]

\newtheorem{Proposition}[Theorem]{Proposition}
\newtheorem{Lemma}[Theorem]{Lemma}
\newtheorem{Corollary}[Theorem]{Corollary}
\newtheorem{Remark}[Theorem]{Remark}
\date{}
\title{Autonomous
Integral Functionals with Discontinuous\\
Nonconvex Integrands:
Lipschitz Regularity of Minimizers,\\
DuBois-Reymond Necessary Conditions,\\
and Hamilton-Jacobi Equations}

\author{Gianni Dal Maso\thanks{ {\it SISSA, via Beirut 2-4,
34014 Trieste, Italy \hfill\break
\null\hspace{2\parindent} e-mail:\/} {\tt dalmaso@sissa.it}}
\hspace{.1cm}
and H\'el\`ene Frankowska\thanks{ {\it CNRS, CREA, Ecole Polytechnique, 1,
Rue Descartes, 75005 Paris, France \hfill\break
\null\hspace{2\parindent} e-mail:\/}
{\tt franko@poly.polytechnique.fr}}
}

\begin{document}

\maketitle

\begin{abstract}
This paper is devoted to the
autonomous Lagrange problem of the calculus of variations with
a discontinuous Lagrangian. We prove that every minimizer is
Lipschitz continuous if the Lagrangian is coercive and
locally bounded. The main difference with respect to
the previous works in the literature is
that we do not assume that
the Lagrangian is convex in the velocity. We also show that, under some
additional assumptions, the DuBois-Reymond necessary condition
still holds in the discontinuous case. Finally, we apply
these results to deduce that the value function of the Bolza
problem is locally Lipschitz and satisfies (in a generalized sense)
a Hamilton-Jacobi equation.
\end{abstract}

\vspace{ 5 mm}
{\small


\vspace{ 5 mm}

}

\newpage
\section{Introduction}\label{intro}

In this paper we study the Lipschitz continuity of the solutions
to the Lagrange problem of the calculus of
variations
\begin{equation}\label{dm300}
\min \left\{ \int_{a}^{b}L(y(t),y'(t))dt \;|\;y \in
W^{1,1}(a,b; \R^n),\;\; y(a)=x_a,\;\; y(b)= x_b\right\},
\end{equation}
where the Lagrangian $L:\R^n \times \R^n \rightarrow \R_+$  is
   a Borel function having a superlinear growth with respect
to the second variable, i.e.,
there exists a function $ \Theta \colon  \R^n \rightarrow \R_+ $,
with
\begin{equation}
   \label{dm302}
\lim_{|u| \rightarrow  \infty } \frac{\Theta (u)}{|u|} = + \infty,
\end{equation}
such that
\begin{equation}
    \label{dm304}
\forall \; (x,u) \in \R^n \times \R^n, \;\; L (x,u) \geq \Theta
(u) .
\end{equation}
We assume also that
$L$ is bounded in a neighborhood of each point of $\R^n \times
\{0\}$, i.e.,
\begin{equation} \label{dm301}
   \forall \; x_0\in\R^n,\; \exists\; M>0,\;\exists \; r>0,\; \forall \;
(x,u) \in B (x_0,r) \times B (0,r),\; L ( x,u) \leq M,
\end{equation}
where $B (x_0,r)$ is the closed ball with center $x_0$ and radius~$r$.

The existence of a solution to (\ref{dm300}) is an easy consequence of
the direct method of the calculus of variations when the functional
\begin{equation}\label{eq150}
{\cal L}_a^b(y ):=
\int_{a}^{b} L(y(t),y'(t))dt
\end{equation}
is  sequentially weakly lower
semicontinuous  on $W^{1,1}(a,b;\R^n)$.

By the classical results of Olech \cite{Olech} and Ioffe
\cite{Ioffe},
a standard assumption for the semicontinuity of ${\cal L}_a^b$ is that
$L$ is lower semicontinuous on $\R^n\times\R^n$ and $L(x,\cdot)$ is
convex on $\R^n$ for every $x\in \R^n$, but these conditions are not
necessary for the lower semicontinuity of ${\cal L}_a^b$ (see, e.g.,
\cite{DG-But-DM}).

Recently Amar, Bellettini and Venturini have
proved in \cite{amar} that any integral functional of the form
(\ref{eq150}), satisfying suitable growth conditions,
has a lower semicontinuous envelope $\overline{\cal L}{}_a^b$
that can be written
as
\begin{equation}\label{calL}
\overline{\cal L}{}_a^b(y)=\int_{a}^{b} L^+(y(t),y'(t))dt\,,
\end{equation}
where $L^+$ is an integrand depending on $L$
(see (\ref{feb30})). If $L$ is continuous, then
$L^+$ coincides with the convexification $L_0$ of $L$ with respect to $u$,
but, if $L$ is discontinuous, one can prove only that the
function $L^+(x,\cdot)$ is convex for a.e.\ $x\in\R^n$, and there are
examples where $L^+(\cdot,u)$ is not lower semicontinuous on $\R^n$.

This shows that there are problems of the form
(\ref{dm300}) which admit a solution even if $L$ is not convex in $u$
(nor lower semicontinuous in $x$), and provides a motivation for the
study of the Lipschitz continuity of the solutions of (\ref{dm300})
without convexity hypotheses.

If $L(x,\cdot)$ is convex for every $x\in\R^n$,
it was proved by Ambrosio, Ascenzi, and Buttazzo in \cite{ambr}
that every minimizer of (\ref{dm300}) is Lipschitz continuous. This kind of
results goes back to Tonelli \cite{t1,t2} for smooth Lagrangians, and
is the first step to prove, under some additional conditions on $L$,
that all minimizers are smooth
(see, e.g., \cite[Section 2.6]{ces83}). Note that, in general,
when the Lagrangian is time
dependent, the problem may have no Lipschitz minimizer (see
\cite{bm85} and \cite{cla-vin}).

The aim of Section~\ref{main} of the present paper
is to show that the convexity hypothesis can be removed
from \cite{ambr}. Assuming only (\ref{dm302}),  (\ref{dm304}), and
(\ref{dm301}), we prove that all minimizers of (\ref{dm300})
are still Lipschitz continuous (Theorem~\ref{dm306}), and provide an
estimate of the Lipschitz constant if, in addition, $L$ is locally bounded
(Theorem~\ref{dm306a}).

  If $L$ is continuous,
then every minimizer $y$ of (\ref{dm300}) is also a minimizer of
the same problem with $L$ replaced by its convexification $L_0$ with
respect to $u$, so that
the Lipschitz continuity of $y$ follows from \cite{ambr}. But, if
$L$ is discontinuous,
we can only say (under suitable growth conditions) that $y$ is a minimizer
of (\ref{dm300}) with $L$ replaced by $L^+$, and we know that
$L^+(x,\cdot)$ is convex only for a.e.\ $x\in\R^n$. For this reason we can
not apply the results of \cite{ambr}. On the other hand,
the proof of \cite{ambr} is
based on an extension
of the DuBois-Reymond
necessary condition, which is not always valid
when $L(x,\cdot)$ is not convex. Therefore we need different arguments.

As in \cite{ambr}, we begin by proving (Lemma~\ref{dm309})
that if $y$ is a minimizer of
(\ref{dm300}), then the function $\psi(t):=t$ is a minimizer
of the problem
\begin{displaymath}
\min \left\{ \int_{a}^{b}f(t,\psi'(t))dt \;|\;\psi \in
W^{1,1}(a,b),\;\; \psi(a)=a, \;\; \psi(b)=b\right\},
\end{displaymath}
where
\begin{equation}
\label{dm307} f (t,v):= \left\{ \begin{array}{lll} L  ( y
(t), {y' (t)}/{v} )v & \mbox{ if}  &
v >
\frac{1}{2},
\\
\noalign{\vskip 3pt}
  + \infty &  \mbox{ if}  &
v \leq \frac{1}{2}.
\end{array} \right.
\end{equation}
Then we show (Lemma~\ref{dm374}) that $\psi(t):=t$ is a minimizer
of the problem
\begin{displaymath}
\min \left\{ \int_{a}^{b}f_0(t,\psi'(t))dt \;|\;\psi \in
W^{1,1}(a,b),\;\; \psi(a)=a, \;\; \psi(b)=b\right\},
\end{displaymath}
where
$f_0= \overline {\rm co}_v f$ is the lower semicontinuous convex
envelope of $f$ with respect to $v$.
This implies (Lemma~\ref{dm376}) that there exists
a constant  $c \in \R$ such that
\begin{displaymath}
d^{\,l}_vf_0 (t,1) \leq c \leq d^{\,r}_v f_0 (t,1) \;\; \mbox{
for a.e.}\;\; t \in [ a,b],
\end{displaymath}
where $d^{\,l}_v$ and $d^{\,r}_v$ denote the left and right
derivatives with respect to $v$.

These inequalities, together with (\ref{dm302}),  (\ref{dm304}), and
(\ref{dm301}), are used to obtain
a bound on the Lipschitz constant of $y$
(Theorem~\ref{dm306}), which is locally
uniform (with respect to the data of the problem) if $L$ is locally bounded
(Theorem~\ref{dm306a}).

In Section~\ref{DuB} we obtain
some extensions of the DuBois-Reymond
necessary condition.
When $L(x,\cdot)$ is not convex
this condition is not always satisfied, and we
propose some additional assumptions on $L$,
which hold true, for instance, when $L(x,\cdot)$ is semiconvex
or differentiable.
Under these assumptions
we show (Theorems~\ref{dm3pr70} and~\ref{dm363})
that, if $y$ is a
minimizer,
then there
exists a constant  $c \in \R$ such that
\begin{displaymath}
   c \in  L (y (t),y' (t))-
   \left\langle   \partial _{u} ^-L (y (t),y' (t)), y' (t)\right\rangle \;\;
\mbox{ for a.e.}\;\; t \in [ a,b],
\end{displaymath}
where $ \partial _{u} ^-L (y (t),y' (t))$ is the subdifferential
of $L(y(t), \cdot)$ at $y'(t)$. More general results of this kind
(Lemma~\ref{dm363n} and Theorem~\ref{dm363x}) are obtained with
different generalized gradients of $L$.

Finally, in Section~\ref{Hamilton} we apply the Lipschitz regularity of
minimizers to study the value function of the Bolza problem:
\begin{equation} \label{value}
   V(t,x):= \inf \left\{\int_{0}^{t} L (y (s), y' (s))ds + \varphi  (y
(t)) \;|\; y\in W^{1,1}(0,t;\R^n),\  y (0)=x \right\} ,
\end{equation}
where $ \varphi \colon  \R^n \to \R_+ \cup \{+ \infty\}$, $\varphi \not\equiv
+ \infty $, and $L$ is locally bounded, not necessarily convex
with respect to the second variable.

Let $\R_+^{\star}:=\{t\in\R \;|\; t>0\}$ and let
$H$ be the Hamiltonian associated
with~$L$, defined by
\begin{equation}\label{H}
   H (x,p):= \sup_{u \in \scriptR^n} \left( \left\langle   p,u\right\rangle
   -L(x,u) \right),
\end{equation}
i.e., $H (x, \cdot )$ is the Legendre-Fenchel transform of $L (x,\cdot )$.

Assuming that for all $(t,x)
\in \R_+^{\star} \times \R^n$ the infimum in (\ref{value}) is attained,
we prove that $V$ is locally Lipschitz on $\R_+^{\star} \times
\R^n$ (Theorem~\ref{liploc})
and solves the Hamilton-Jacobi equation
\begin{equation}\label{HJ}
V_t + H(x,-V_x)=0
\end{equation}
in a generalized sense (Theorem~\ref{thm2}).
When
$\varphi$  is lower semicontinuous, we also provide a comparison
result for lower semicontinuous
subsolutions of~(\ref{HJ}), which
characterizes the
value function as the maximal lower semicontinuous subsolution
of~(\ref{HJ}) (Theorem~\ref{thm1}).

We conclude the paper with two results
(Theorems~\ref{dm391} and~\ref{dm391b}) which show the relationships
between minimizers of (\ref{value}) and contingent derivatives of the
value function.

\section{Lipschitz Regularity of Minimizers}\label{main}

Let $L\colon \R^n \times \R^n \rightarrow \R_+$ be a Borel function,
let  $[a,b]$ be a bounded closed interval in $\R$, and let  $y \in W^{1,1}
(a,b; \R^n)$ be a function such that
\begin{equation}
   \label{dm305}
\int_{a}^{b}L (y (t), y' (t))dt \leq \inf_{z \in S(y)}
\int_{a}^{b} L (z (t), z' (t))dt  < + \infty,
\end{equation}
where
$
S(y):=\{z \in W^{1,1} (a,b; \R^n) \; | \; z (a)= y (a),
\;\;z (b)= y (b)\}$.

The main results of this section are the following two theorems.

\begin{Theorem} \label{dm306}
Let $L\colon \R^n \times \R^n \rightarrow \R_+$ be a Borel function
which satisfies   (\ref{dm302}),  (\ref{dm304}), (\ref{dm301}), and
let  $y \in W^{1,1} (a,b; \R^n)$ be a function which satisfies
(\ref{dm305}). Then $y$ is Lipschitz continuous.
\end{Theorem}

When $L$ is locally bounded on $\R^n \times \R^n $, we obtain a
uniform estimate of the Lipschitz constant of every minimizer.

\begin{Theorem} \label{dm306a}
Let $\Theta\colon  \R^n \to \R_+$ be a function satisfying
(\ref{dm302}), let
$\Psi \colon  \R_+ \to \R_+$ be a nondecreasing function, and let
$A$, $B$, $\alpha$, $\beta >0$. Then there exists a constant $K
=K(\Theta, \Psi, A,B, \alpha,\beta)$ with the
following property: if $L\colon \R^n \times \R^n \rightarrow \R_+$ is
any Borel function
satisfying (\ref{dm304}) and
\begin{equation}
\label{dm379}
\forall \; R>0, \quad \sup_{|x| \leq R, \,|u|\le R} L(x,u) \leq \Psi(R),
\end{equation}
and $y \in W^{1,1} (a,b; \R^n)$ satisfies (\ref{dm305}) and
\begin{eqnarray}
   \label{dm369a}
   &\displaystyle\vphantom{\int_{a}^{b}}
\inf_{a \leq t \leq b} |y(t)| \leq A ,
\\
   \label{dm370a}
&\displaystyle
\int_{a}^{b} L(y(t),y'(t))dt \leq B ,
\\
   \label{dm371}
   &\displaystyle\vphantom{\int_{a}^{b}}
\alpha \leq  b-a \leq \beta ,
\end{eqnarray}
then $y$ is Lipschitz continuous with Lipschitz constant bounded by~$K$.
\end{Theorem}

To prove Theorems \ref{dm306} and \ref{dm306a}
we need some technical lemmas.

Let us fix a function $y\in W^{1,1}(a,b;\R^n)$
which satisfies (\ref{dm305}).
As in  \cite{ambr} we use the auxiliary function $f \colon [a,b]
\times \R \to  [0, + \infty ]$ defined by (\ref{dm307}),
which turns out to be  $ {\cal L}_1 \times {\cal B}_1$-measurable,
where $ {\cal L}_1  $  and $ {\cal B}_1$ denote the
$\sigma$-algebras of Lebesgue measurable subsets of  $[a,b]$ and of
Borel subsets of  $\R$, respectively. {}From
(\ref{dm307})  and (\ref{dm305}) it follows that
\begin{equation}
   \label{dm308}
\int_{a}^{b} f (t,1) dt < + \infty .
\end{equation}

The following lemma is well known
(see, e.g., \cite[p.~46]{ces83}). We write the proof only to give a
self-contained presentation of the arguments used to obtain
Theorems \ref{dm306} and \ref{dm306a}.
\begin{Lemma}
\label{dm309}
We have
\begin{equation}
   \label{dm310}
\int_{a}^{b}f (t,1)dt \leq \int_{a}^{b}f (t, \psi ' (t))dt
\end{equation}
for every  $\psi \in W^{1,1} (a,b)$ such that  $ \psi (a)=a$ and  $\psi
(b)=b$.
\end{Lemma}

\noindent
{\bf Proof} --- \hspace{ 2 mm}
Let us fix $\psi \in W^{1,1}
(a,b)$, with  $ \psi (a)=a$ and  $\psi (b)=b$, such that the right
hand side of   (\ref{dm310}) is finite. Then $\psi ' (t) >
\frac{1}{2}$ for almost all  $ t \in [ a,b]$. Thus  $ \psi $ is
increasing and  $| \psi (t)- \psi  (s)| \geq \frac{1}{2} |t-s|$
for all $s,t \in [ a,b]$. Therefore the inverse function $
\psi^{-1} \colon  [ a,b] \to  [ a,b]$ is Lipschitz continuous with Lipschitz
constant 2. These properties imply that
\begin{equation}
   \label{dm311}
N \in {\cal L}_1, \;\; |N|=0  \Longleftrightarrow \psi (N) \in
{\cal L}_1,\;\; | \psi  (N)| =0 ,
\end{equation}
where  $|\cdot |$ denotes the Lebesgue measure, and that  $z\colon = y
\circ \psi ^{-1}$ belongs to $W^{1,1} (a,b;\R^n)$ and satisfies $z
(a)=y (a)$ and  $z (b)= y (b)$. Using  (\ref{dm311}) and the chain
rule one proves easily that $z' (t)=y' ( \psi ^{-1} (t))/ \psi ' (
\psi ^{-1} (t))$ for a.e.\ $t \in [ a,b]$. Thus, after the change
of variables  $s= \psi ^{-1}  (t)$, one gets
\begin{displaymath}
   \int_{a}^{b} L (z (t),z' (t))dt= \int_{a}^{b}L ( y (s), {y'
(s)}/{ \psi ' (s)} ) \psi ' (s) ds ,
\end{displaymath}
which, together with  (\ref{dm305}) and  (\ref{dm307}), yields
(\ref{dm310}).
$\; \; \Box$

\vspace{ 5 mm}

If  $g\colon  \R \rightarrow  [ 0, + \infty ]$ is an arbitrary function, its
lower semicontinuous convex envelope $ \overline {\rm co} \, g \colon  \R
\rightarrow  [ 0, + \infty ]$ is, by definition,  the greatest lower
semicontinuous convex function which is less than or equal to  $g$. It is
well known that the epigraph of  $\overline {\rm co} \, g$ is the closed convex
hull of the epigraph of  $g$, and that, if  $ \overline {\rm co} \,g$ is finite
in a neighborhood of some point  $v \in \R$, then
\begin{equation}
   \label{dm312}
\overline {\rm co}\, g (v)= \inf \left\{ \lambda g(v_1) + (1-
\lambda ) g (v_2)\;|\;\;  ( \lambda ,v_1, v_2) \in A (v) \right\}
,
\end{equation}
where  $A (v)$ is the set of all triples $ ( \lambda ,v_1, v_2) \in \R^3$
with $ 0 < \lambda  <1,\; g (v_1) < + \infty , \; g (v_2) < + \infty ,$ and
$v= \lambda v_1+ (1- \lambda )v_2$.

Let us return to the function  $f$ defined in  (\ref{dm307}), and
let $f_0= \overline {\rm co}_v f$ be its lower semicontinuous
convex envelope with respect to $v$. We observe that for every  $t
\in [a,b]$, the function $v \mapsto  f_0 (t,v)$ is continuous on
$( \frac{1}{2}, + \infty )$, since it is convex and finite on this
set.

Let us prove that for every $v \in ( \frac{1}{2}, + \infty)$ the
function $t \mapsto f_0 (t,v)$ is Lebesgue measurable. Given  $
\gamma \in \R$, by  (\ref{dm312}) the set $\{t \in [ a,b] \;| \;
f_0 (t,v) < \gamma \}$ is the projection onto  $ [ a,b]$ of the
set of all points $ (t, \lambda , v_1,v_2) \in [ a,b] \times (0,1)
\times ( \frac{1}{2}, + \infty ) \times ( \frac{1}{2}, + \infty )$
such that $v = \lambda v_1 + (1 - \lambda )v_2$ and $ \lambda f
(t,v_1) + (1- \lambda ) f (t,v_2) < \gamma $. As this set belongs
to the $ \sigma$-algebra  $ {\cal L} _1 \times {\cal B}_1\times
{\cal B}_1\times {\cal B}_1$, from the projection theorem (see,
e.g., \cite[Theorem 8.3.2]{af90sva}) we conclude that $ \{t \in [
a,b]\;|\; f_0 (t,v) < \gamma \}$ is Lebesgue measurable. This
proves that $t \to f_0 (t,v)$ is Lebesgue measurable, and hence
$f_0$ is a Carath\'{e}odory function on  $ [ a,b] \times (
\frac{1}{2}, + \infty )$.\

\vspace{ 5 mm}
The following lemma is usually proved when $f$ is continuous in $v$, or
satisfies some growth condition.  We give here a detailed proof
to show that we do not need any additional hypothesis.

\begin{Lemma}\label{dm374}
We have
\begin{equation}
\label{dm313}
\int_{a}^{b}f_0(t,1)dt \leq
   \int_{a}^{b}f (t,1)dt \leq \int_{a}^{b} f_0 (t, \varphi  (t))dt
\end{equation}
for every  $ \varphi \in L^1 (a,b)$  with
\begin{equation}
\label{dm314}
   \int_{a}^{b} \varphi  (t) dt =b-a .
\end{equation}
In particular,
\begin{equation}
    \label{dm326}
f (t,1)=f_0 (t,1) \;\; \mbox{for a.e.\ }t \in  [ a,b].
\end{equation}
\end{Lemma}

\noindent
{\bf Proof} --- \hspace{ 2 mm} The first inequality in (\ref{dm313})
follows from the fact that $f_0(t,1)\le f(t,1)$ for every
$t\in[a,b]$. To prove the second inequality
we argue by contradiction. Assume that there exists $ \varphi  \in L^1
(a,b)$, satisfying  (\ref{dm314}), such that
\begin{equation} \label{dm315}
   \int_{a}^{b} f_0 (t, \varphi  (t))dt < \int_{a}^{b} f (t,1)dt .
\end{equation}
As  $f_0 (t,v) = + \infty $ for  $v < \frac{1}{2}$, from  (\ref{dm308})
and  (\ref{dm315}) we obtain that $ \varphi  (t) \geq \frac{1}{2}$ for
a.e.\ $t \in  [ a,b]$. Changing, if needed, $ \varphi $ on a set of measure
zero, we may assume that this inequality is satisfied for every  $t \in  [
a,b]$. If we replace  $ \varphi  (t)$ by  $ \frac{1}{2}  \varphi  (t) +
\frac{1}{2}$, we obtain a new function, still denoted by  $ \varphi (t)$,
which continues to fulfill  (\ref{dm314}) and  (\ref{dm315}) (by
convexity), and, in addition, satisfies the improved inequality $ \varphi
(t) \geq \frac{3}{4}$ for every  $t \in  [ a,b]$. As $f_0$
is a Carath\'{e}odory
function on  $ [ a,b] \times ( \frac{1}{2}, + \infty )$, the function
$t\mapsto f_0(t,\varphi(t))$ is measurable.

Let us fix  $ \varepsilon >0$ such that
\begin{equation}
   \label{dm316}
\int_{a}^{b} [ f_0 (t, \varphi  (t)) + \varepsilon ]dt \; <\;
\int_{a}^{b} f (t,1) dt .
\end{equation}
For every  $t \in  [ a,b]$ we define $ A_{ \varepsilon } (t)$ as the set of
all triples  $ ( \lambda ,v_1,v_2) \in \R^3$ such that  $0 < \lambda <1$,
$v_1 > \frac{1}{2}$, $v_2 > \frac{1}{2}$,
$\lambda v_1+ (1- \lambda )v_2 =
\varphi  (t)$, and
\begin{displaymath}
   \lambda f (t,v_1) + (1- \lambda )f (t,v_2) < f_0 (t, \varphi  (t)) +
\varepsilon .
\end{displaymath}
By  (\ref{dm312}) this set is nonempty for every  $t \in  [ a,b]$.

{}From the measurability properties of  $f$ and
$t\mapsto f_0(t,\varphi(t))$ we deduce
that the graph of the set-valued map $t \leadsto A_{ \varepsilon }
(t)$ belongs to $ {\cal L} _1 \times {\cal B}_1\times {\cal
B}_1\times {\cal B}_1$. By the Aumann--von Neumann selection
theorem (see, e.g., \cite[Theorem III.22]{cava}) there exist
Lebesgue measurable functions $ \mu , \varphi _1, \varphi _2$ on
$ [ a,b]$ such that $ ( \mu (t), \varphi _1 (t), \varphi _2 (t))
\in A_{ \varepsilon } (t)$ for every  $t \in [ a,b]$. {}From the
definition of  $A_{ \varepsilon } (t)$ and from  (\ref{dm316}) we
deduce that
\begin{eqnarray}
   \label{dm317}
&\displaystyle\vphantom{\int_{a}^{b}}
\forall \; t \in [ a,b],\;\; \; 0 < \mu  (t) <1,\; \; \; \varphi _1 (t) >
\frac12,\; \; \;\varphi _2 (t) >
\frac12,
\\
   \label{dm318}
&\displaystyle\vphantom{\int_{a}^{b}}
\forall \; t \in [ a,b],\;\; \; \mu  (t) \varphi _1 (t)+ (1 - \mu
(t)) \varphi _2 (t) = \varphi (t) ,
\\
   \label{dm319}
&\displaystyle\int_{a}^{b} [ \mu  (t) f (t, \varphi _1 (t)) + (1- \mu
(t)) f (t,
\varphi _2 (t))] dt < \int_{a}^{b}f (t,1)dt .
\end{eqnarray}
For every  $k \geq 2$ let
\begin{displaymath}
   E^{ (k)}:= \left\{ t \in  [ a,b]\;|\; \mu (t) \in ( \frac{1}{k+1},
\frac{1}{k}] \cup [ \frac{k-1}{k}, \frac{k}{k+1}) \right\} .
\end{displaymath}
As $0< \mu (t) <1$, the interval  $[ a,b]$ is the union of the sets
$E^{ (k)}$, which are pairwise disjoint. As  $\varphi\in L^1(a,b)$
and $f(\cdot,1)\in L^1(a,b)$ by (\ref{dm308}),
from
(\ref{dm318}) and (\ref{dm319}) we obtain that for every  $k \geq
2$ the functions $ \varphi _i (t)$ and $f (t, \varphi _i (t))$
belong to $L^{1} (E^{ (k)})$.

By the Lyapunov theorem (see, e.g.,  \cite[Theorem 8.6.3 and
Proposition 8.6.2]{af90sva}) there exist two disjoint measurable
sets $ E^{ (k)}_{1}$ and  $ E^{ (k)}_{2}$, with  $E^{ (k)}_{1}
\cup E^{ (k)}_{2}=E^{ (k)}$, such that
\begin{eqnarray}
    \label{dm321}
    &\displaystyle
\sum_{i=1}^{2} \int_{E^{ (k)}_{i}}^{} \varphi _i (t) dt = \int_{E^{
(k)}}^{} [  \mu  (t) \varphi _1 (t) + (1- \mu  (t)) \varphi _2 (t)] dt,
\\
    \label{dm322}
    &\displaystyle
\sum_{i=1}^{2} \int_{E^{ (k)}_{i}}^{} f(t, \varphi _i (t)) dt =
\int_{E^{ (k)}}^{} [  \mu  (t) f(t,\varphi _1 (t)) + (1- \mu
(t))f(t, \varphi _2 (t))] dt .
\end{eqnarray}
Let $\displaystyle
   E_1:= \bigcup_{k=2}^{ \infty }E^{ (k)}_{1}$ and $
   \displaystyle E_2:=
   \bigcup_{k=2}^{\infty }E^{ (k)}_{2}$.
By  (\ref{dm314}), (\ref{dm318}),   (\ref{dm319}), (\ref{dm321}),
(\ref{dm322}) we obtain
\begin{eqnarray}
   \label{dm323}
&  \displaystyle
\int_{E_1}^{} \varphi _1 (t)dt + \int_{E_2}^{} \varphi _2 (t) dt=
\int_{a}^{b} \varphi  (t)dt =b-a,
\\
    \label{dm324}
  &    \displaystyle
\int_{E_1}^{} f (t, \varphi _1 (t))dt + \int_{E_2}^{}f (t, \varphi
_2 (t)) dt < \int_{a}^{b} f (t,1)dt .
\end{eqnarray}
Let $ \varphi_3 \in L^{1} (a,b)$ be the function defined by $ \varphi
_3:= \varphi_1$ on  $E_1$ and  $ \varphi_3 := \varphi _2$ on  $E_2$, and
let  $ \psi $ be the primitive of  $ \varphi_3$ with  $ \psi  (a)=a$. By
(\ref{dm323}) we have also  $ \psi  (b)=b$, while (\ref{dm324}) gives
\begin{displaymath}
   \int_{a}^{b}f (t, \psi ' (t))dt < \int_{a}^{b} f (t,1)dt,
\end{displaymath}
which contradicts   (\ref{dm310}) and concludes the proof of
(\ref{dm313}).

As  $f_0 \leq f$ and $f(\cdot,1)\in L^1(a,b)$ by (\ref{dm308}),
if we take $ \varphi  \equiv 1$ in  (\ref{dm313}) we
get~(\ref{dm326}).
$\; \; \Box$

\begin{Remark}{\rm As  $v\mapsto f_0 (t,v)$ is convex and finite,
     for every $t \in  [
a,b]$ there exist the limits
\begin{eqnarray}
\label{dm327}
&\displaystyle
d^{\,l}_v f_0 (t,1):= \lim_{v \rightarrow 1-}  \frac{f_0 (t,v)-f_0
(t,1)}{v-1} = \sup_{v<1} \frac{f_0 (t,v)-f_0 (t,1)}{v-1},
\\
\label{dm328}
&\displaystyle
d^{\,r}_v f_0 (t,1):= \lim_{v \rightarrow 1+}  \frac{f_0 (t,v)-f_0
(t,1)}{v-1} = \inf_{v>1} \frac{f_0 (t,v)-f_0 (t,1)}{v-1},
\end{eqnarray}
and we have $ - \infty < d^{\,l}_v f_0 (t,1) \leq d^{\,r}_v f_0 (t,1) < +
\infty $.
}
\end{Remark}

For the sake of completeness, we give now a new elementary proof of a
particular case of Theorem 3.1 of  \cite{ambr}.
\begin{Lemma}\label{dm376}
There exists a constant  $c \in \R$ such that
\begin{equation}
   \label{dm329}
d^{\,l}_vf_0 (t,1) \leq c \leq d^{\,r}_v f_0 (t,1) \;\; \mbox{
for a.e.}\;\; t \in [ a,b].
\end{equation}
\end{Lemma}

\noindent
{\bf Proof} --- \hspace{ 2 mm}
We argue by contradiction. If (\ref{dm329}) does not hold, then there
exists $ \alpha \in \R$ such that
\begin{displaymath}
   \mathop{\rm ess \, sup }_{t \in  [ a,b]}d^{\,l}_v f_0 (t,1) > \alpha >
\mathop{\rm ess \, inf } _{t \in
[ a,b]}d^{\,r}_v f_0 (t,1) .
\end{displaymath}
Then the sets
\begin{displaymath}
   A_l := \left\{ t \in [ a,b] \;| \; d^{\,l}_v f_0 (t,1) > \alpha  \right\}
\;\; \& \;\;  A_r :=\left\{ t \in [ a,b] \;| \; d^{\,r}_v f_0 (t,1) <
\alpha  \right\}
\end{displaymath}
are disjoint and have positive measure. By  (\ref{dm327}) for every  $t \in
A_l$ the set
\begin{displaymath}
   E_l (t):= \left\{ v \in (  \frac{1}{2},1 ) \;|\; \frac{f_0
(t,v) -f_0 (t,1)}{v-1} > \alpha  \right\}
\end{displaymath}
is nonempty. Since the graph of the set-valued map $t \leadsto E_l
(t)$ belongs to $ {\cal L}_1 \times {\cal B}_1$, by the Aumann--von
Neumann selection theorem (see, e.g.,  \cite[Theorem III.22]{cava})
there exists a measurable function  $ \delta _l \colon  A_l \to (0,
\frac{1}{2})$ such that
\begin{equation}
    \label{dm330}
\forall \; t \in A_l,\;\; f_0 (t,1- \delta _l (t)) -f_0 (t,1) < -
\alpha \delta _l (t).
\end{equation}
Similarly, using  (\ref{dm328}) we can prove that there exists a
measurable function  $ \delta _r \colon  A_r \to (0, \frac{1}{2})$
such that
\begin{equation}
    \label{dm331}
\forall \; t \in A_r,\;\; f_0 (t,1+ \delta _r (t))-f_0 (t,1) <
\alpha \delta _r (t).
\end{equation}
Let us define
\begin{equation}
    \label{dm332}
c_l:= \left[  \int_{A_l}^{} \delta _l (t) dt\right]^{-1}\;\;\&\;\;
c_r:= \left[  \int_{A_r}^{} \delta _r (t) dt\right]^{-1},
\end{equation}
and let  $ \varphi  (t):=-c_l \delta _l (t)$ for $t \in A_l,\; \varphi
(t)=c_r \delta _r (t)$ for $t \in A_r$, and  $ \varphi  (t)=0$ otherwise.
Then by  (\ref{dm332}) we have  $ \int_{a}^{b} \varphi  (t) dt=0$, and,
by   (\ref{dm313}), for every  $ \varepsilon >0$ this implies
\begin{displaymath}
   \int_{a}^{b} \left[  f_0 (t,1+ \varepsilon  \varphi  (t)) - f_0
(t,1)\right] dt \geq 0 ,
\end{displaymath}
which is equivalent to
\begin{equation}
    \label{dm333}
\int_{A_l}^{} \left[  f_0 (t,1- \varepsilon c_l \delta _l (t)) -
f_0 (t,1)\right] dt + \int_{A_r}^{} \left[  f_0 (t,1+ \varepsilon
c_r \delta _r (t)) - f_0 (t,1)\right] dt \geq 0 .
\end{equation}
By the monotonicity property of the difference quotient of a convex
function, using   (\ref{dm330}) we obtain for  $ \varepsilon c_l <1$
\begin{equation}   \label{dm334}
f_0 (t,1- \varepsilon c_l \delta _l (t))
-f_0 (t,1) \leq  \varepsilon c_l \left[  f_0 (t, 1- \delta _l (t))
- f_0 (t,1)\right]  < - \alpha  \varepsilon c_l \delta _l (t)
\end{equation}
for every $t \in A_l$.
Similarly,  for $ \varepsilon c_r < 1$ we obtain, using (\ref{dm331}),
\begin{equation} \label{dm3e33}
f_0 (t,1+ \varepsilon c_r \delta _r (t))-f_0
(t,1) \leq \varepsilon c_r \left[ f_0 (t,1+  \delta _r (t))-f_0 (t,1)
   \right]
   < \alpha  \varepsilon c_r \delta _r (t)
\end{equation}
for every $t\in A_r$.
{}From (\ref{dm332})--(\ref{dm3e33}) it follows that
\begin{displaymath}
   0 < - \alpha  \varepsilon c_l \int_{A_l}^{} \delta _l (t) dt + \alpha
\varepsilon c_r \int_{A_r}^{} \delta _r (t) dt=0.
\end{displaymath}
This contradiction proves  (\ref{dm329}).
$\; \; \Box$

\vspace{5 mm}
\noindent
{\bf Proof of Theorem \ref{dm306}} --- \hspace{ 2 mm} By
(\ref{dm328}) and (\ref{dm329}) there exists a constant $c\in \R$ such that
for a.e.\ $t \in [ a,b]$ and every
$ \varepsilon  \in (0,1)$ we have
\begin{displaymath}
   c \leq d^{\, r}_v f_0 (t,1) \leq \frac{f_0 (t, 2- \varepsilon ) -f_0
(t,1)}{1- \varepsilon },
\end{displaymath}
hence
$
   (1- \varepsilon )c \leq f_0 (t,2- \varepsilon ) -f_0 (t,1)
$,
which implies
\begin{equation}
   \label{dm3e34}
(1- \varepsilon )c + \varepsilon f_0 (t,1) \leq f_0 (t,2-
\varepsilon ) - (1- \varepsilon )f_0 (t,1) .
\end{equation}
By convexity we have
\begin{equation}
   \label{dm3e35}
f_0 (t,2- \varepsilon ) \leq  \varepsilon f_0 ( t, {1}/{
\varepsilon } ) + (1- \varepsilon )f_0 (t,1),
\end{equation}
so that we obtain from  (\ref{dm3e34}) and  (\ref{dm3e35})
\begin{equation}
   \label{dm3e36}
(1- \varepsilon ) c + \varepsilon f_0 (t,1) \leq \varepsilon f_0
( t, {1}/{ \varepsilon } ) \leq  \varepsilon f (t,
{1}/{ \varepsilon }) .
\end{equation}
By  (\ref{dm304}) and  (\ref{dm307}) for every  $v>0$ we have
\begin{equation}
   \label{dm3e37}
f (t,v) \geq L  ( y (t), {y' (t)}/{v} )v \geq
\Theta ( {y' (t)}/{v} )v \geq \overline{\rm co}\,
\Theta ( {y' (t)}/{v} )v,
\end{equation}
where $\overline{\rm co}\, \Theta $ is the lower semicontinuous convex
envelope of  $\Theta$, which still satisfies (\ref{dm302}). Since
the function $v \mapsto \overline{\rm co}\, \Theta ( {y' (t)}/{v}
)v$ is convex for $v>0$, from (\ref{dm3e37}) we deduce that
\begin{equation}
   \label{dm3e38}
\forall \; v>0,\;\; f_0 (t,v) \geq \overline{\rm co}\, \Theta(
{y' (t)}/{v} )v .
\end{equation}
{}From  (\ref{dm307}),  (\ref{dm3e36}), and  (\ref{dm3e38}) we
obtain
\begin{equation}
   \label{dm3e39}
(1- \varepsilon )c + \varepsilon \, \overline{\rm co}\, \Theta (
y' (t) ) \leq L (y (t), \varepsilon y' (t))
\end{equation}
for a.e.\ $t \in [ a,b]$ and every  $ \varepsilon  \in (0,1)$.

Let us now fix  $t \in [ a,b]$ such that  (\ref{dm3e39}) holds and
$|y' (t)| >2$. For $ \nu \in (0,1]$ let  $ \varepsilon  (t) =
\frac{ \nu }{|y' (t)|} < \frac{1}{2}$. By  (\ref{dm3e39}) we have
\begin{equation}
   \label{dm3e40}
\min \{c,0\} + \nu \frac{\overline{\rm co}\, \Theta ( y' (t) )}{|y'
(t)|} \leq L \left( y (t), \nu \frac{y' (t)}{|y' (t)|} \right).
\end{equation}
Since $y \in W^{1,1} (a,b; \R^n)$, there exists $R>0$ such that
\begin{equation}
   \label{dm3e41}
\forall \; t \in [ a,b],\;\;|y (t)| \leq R.
\end{equation}
Since $B(0,R)$ is compact, from  (\ref{dm301}) we know that
\begin{equation}
   \label{dm3e42}
\exists \; M>0,\; \exists \; r \in (0,1],\;\; \forall \; u \in B
(0,r),\;\; L (y(t),u) \leq M .
\end{equation}
Choosing $ \nu =r$, from  (\ref{dm3e40}),  (\ref{dm3e41}), and
(\ref{dm3e42}) we get
\begin{equation}
   \label{dm3e43}
\min \{c,0\} + r \frac{\overline{\rm co}\, \Theta \left( y' (t) \right)}{|y'
(t)|} \leq M.
\end{equation}
Since $\overline{\rm co}\, \Theta$ satisfies   (\ref{dm302}), by
(\ref{dm3e43}) there exists a constant $C =C(\Theta, c, r,M) \geq 2$,
depending only on $\overline{\rm co}\, \Theta$,
$c$, $r$,  and $M$, such that for a.e.\ $t \in [ a,b]$ with $|y' (t)| >2$
we have
\begin{equation}
   \label{dm3e44}
|y' (t)| \leq C.
\end{equation}
As $C \geq 2$, inequality (\ref{dm3e44}) holds also when  $|y'
(t)| \leq 2$. $\; \; \Box$

\vspace{5 mm}
\noindent
{\bf Proof of Theorem \ref{dm306a}} --- \hspace{ 2
mm} By Lemma \ref{dm376} there exists $c$ such that
\begin{equation}
   \label{dm377}
d^{\,l}_vf_0 (t,1) \leq c \leq d^{\,r}_v f_0 (t,1) \;\; \mbox{
for a.e.}\;\; t \in [a,b] .
\end{equation}
By (\ref{dm304}) and
(\ref{dm370a})  we have
\begin{equation}\label{123}
\int_{a}^{b} \Theta (y'(t)) dt \leq B ,
\end{equation}
which, by (\ref{dm302}), gives
\begin{displaymath}
\int_{a}^{b} |y'(t)| dt\le M_1,
\end{displaymath}
for a constant $M_1=M_1(\Theta, B, \beta)>0$. This inequality,
together with (\ref{dm369a}), yields $|y (t)| \leq R $ for every
$t\in[a,b]$,
with $R=R(\Theta, A,B,\beta)=A+M_1$.

We next provide an estimate of $c$ from below.
{}From (\ref{123}) and (\ref{dm371}) it follows that
\begin{equation}
   \label{dm378}
\alpha \; \mathop{\rm ess\,inf}_{t \in [a,b]} \Theta (y'(t)) \leq
\int_{a}^{b} \Theta (y'(t)) dt \leq B .
\end{equation}
By
(\ref{dm302}) there exists $M_2=M_2(\Theta, \alpha, B)$ such that for a set
$\Omega_y \subset [a,b]$ of positive measure
\begin{equation}
   \label{dm380}
\forall \; t \in \Omega_y,\;\; |y'(t)| \leq M_2 .
\end{equation}
This implies that
\begin{equation}
   \label{dm381}
\forall \; t \in \Omega_y,\;\; L(y(t), \textstyle \frac{4}{3}y'(t))
\leq \Psi (R+ 2 M_2) .
\end{equation}

Since by (\ref{dm327})
\begin{displaymath}
c \geq \sup_{v <1}  \frac{f_0(t,v)
-f_0(t,1)}{v-1}
\end{displaymath}
and since $f_0(t,1) \geq 0$ for almost all $t
\in [a,b]$, we get, setting $v= \frac{3}{4}$,
\begin{equation}
   \label{dm375}\textstyle
c \geq -4f(t,  \frac{3}{4})= -3L(y(t),
\frac{4}{3} y'(t)) \geq -3\Psi (R+ 2 M_2) .
\end{equation}

We now return to the last part of the proof of Theorem \ref{dm306}
with $\nu=r=1$ and $M=\Psi(R+1)$. As the constant $C$ which
appears in (\ref{dm3e44}) depends on $c$ in a decreasing way, and
$c\ge -3\Psi (R+ 2 M_2)$, it is enough to set $K=C(\Theta, -3\Psi
(R+ 2 M_2), 1, \Psi(R+1))$.
$\;\;\Box$

\section{DuBois-Reymond  Necessary Conditions}\label{DuB}

Let $L\colon \R^n \times \R^n \rightarrow \R_+$ be a Borel function,
let $y \in W^{1,1} (a,b; \R^n)$ be a function which satisfies
(\ref{dm305}), and let
$g \colon  [ a,b] \times \R \rightarrow [0, + \infty ]$ be the function
defined by
\begin{equation}
   \label{dm3e45}
g ( t,v):=  \left\{ \begin{array}{lll} L (y (t), vy' (t)) & \mbox{
if\ \ }  0<v<2,\\ + \infty  & \mbox{  otherwise.} &
\end{array} \right.
\end{equation}
By  (\ref{dm307}) we have $f (t,v)=g (t, \frac{1}{v})v$ for every
$v \ne 0$. Let $g_0 := \overline {\rm co}_v g $ be the lower
semicontinuous convex envelope of  $g$ with respect to  $v$. As
the functions $v \mapsto f_0 (t, v)$ and  $v \mapsto g_0 (t,
\frac{1}{v})v$ are lower semicontinuous and convex for $v>0$, we
deduce that
\begin{displaymath}
   \forall \;  t \in [ a,b], \; v>0,\;\; f_0 (t,v)=g_0 \left( t,
{1}/{v} \right)v .
\end{displaymath}
Therefore we obtain from  (\ref{dm326})
\begin{equation}
   \label{dm3e46}
g_0 (t,1)=g (t,1)=L (y (t), y' (t)) \;\; \mbox{ for a.e.}\;\; t \in [ a,b].
\end{equation}
Furthermore $t \mapsto g_0 (t ,v)$ is measurable. To prove this fact
it is enough to adapt the arguments used for $f_0$
in the proof of Lemma \ref{dm309}.

Let us define $d^{\, l}_vg_0 (t,1)$ and  $d^{\,r}_vg_0 (t,1)$ as in
(\ref{dm327}) and (\ref{dm328}). It is easy to prove that
\begin{displaymath}
   d^{\,l}_v f_0 (t,1)=g_0 (t,1)- d^{\,r}_vg_0 (t,1)\;\;\&\;\;
   d^{\,r}_v f_0 (t,1)=g_0 (t,1)- d^{\,l}_vg_0 (t,1).
\end{displaymath}
Therefore, by  (\ref{dm329}) there exists a constant $c \in \R$ such that
\begin{equation}
   \label{dm3e47}
d^{\, l}_v g_0 (t,1)\le L (y (t), y' (t))-c\le d^{\,r}_v g_0 (t,1)
\;\; \mbox{  for a.e.\ }t \in [ a,b].
\end{equation}
Notice that if  $u \mapsto L (y (t),u)$ is differentiable at $y' (t)$,
with gradient $ \nabla _u L (y (t),y' (t))$, then by (\ref{dm3e45})
the function $v \mapsto g (t,v)$ has a derivative at  $v=1$ which is
equal to
$ \left\langle    \nabla _u L (y (t),y' (t)),y'(t)\right\rangle $.
By (\ref{dm3e46}) this implies
\begin{equation}
   \label{dm3e48}
d^{\,l}_vg_0 (t,1)=d^{\, r}_v g_0 (t,1)= \left\langle   \nabla _uL (y
(t),y' (t)),y' (t)\right\rangle,
\end{equation}
and from  (\ref{dm3e47}) and (\ref{dm3e48}) we obtain the
DuBois-Reymond  necessary condition
\begin{displaymath}
L (y (t),y' (t))- \left\langle   \nabla _u L (y (t),y' (t)),y'
(t)\right\rangle =c\;\; \mbox{for a.e.}\;\; t \in [ a,b] .
\end{displaymath}

Our aim is to derive similar results when  $L (y (t), \cdot )$ is not
differentiable.
All our extensions of the
DuBois-Reymond  necessary condition
(Theorems~\ref{dm3pr70}, \ref{dm363},
and~\ref{dm363x})
are based on the following lemma.

\begin{Lemma}\label{dm363n}
Let $L\colon \R^n \times \R^n \rightarrow \R_+$ be a Borel function,
let  $y \in W^{1,1} (a,b; \R^n)$ be a function which satisfies
(\ref{dm305}), and let
$\psi\colon [a,b] \times \R^n \to
\R$ be a Carath\'eodory function, with $\xi\mapsto \psi(t,\xi)$
convex and positively
homogeneous of degree one, such that for a.e.\ $t \in [ a,b]$
\begin{equation}
\label{dm360'n}
  -d^{\,l}_vg_0 (t,1)\le \psi(t, -y' (t)) \quad\&\quad
d^{\,r}_vg_0 (t,1)\le \psi(t, y' (t)).
  \end{equation}
Then there exist a constant  $c \in \R$ and a measurable function
$p \colon  [ a,b] \rightarrow \R^n$ such that for a.e.\ $t \in [ a,b]$
\begin{eqnarray}
   \label{dm362*}
& p (t) \in \partial_\xi \psi(t,0) ,
\\
\label{dm3a60*}
& L (y (t),y' (t))- \left\langle   p (t), y' (t)\right\rangle = c ,
\end{eqnarray}
where $\partial_\xi \psi(t,0)$ denotes the subdifferential of the convex
function $\psi(t,\cdot)$ at~$0$.
\end{Lemma}

\noindent
{\bf Proof} --- \hspace{ 2 mm} Since $\psi(t, \cdot)$ is convex and
positively homogeneous of degree one,
\begin{equation}\label{dm3650}
\forall \; \xi \in \R^n,\;\; \max_{q \in \partial _{\xi} \psi(t,0)}
\langle q,\xi\rangle = \psi(t,\xi).
\end{equation}
  By  (\ref{dm3e47}) and (\ref{dm360'n})
  there exists a constant  $c$ such that for a.e.\ $t \in [ a,b]$
\begin{displaymath}
  -\psi(t,-y'(y))\leq L (y (t),y' (t))-c \leq \psi(t,y'(t)),
\end{displaymath}
so that  (\ref{dm3650}) implies
\begin{displaymath}
  \min _{q \in \partial _{\xi} \psi(t,0) } \langle q,y'(t)\rangle
  \leq L (y (t),y' (t)) - c  \leq
  \max _{q\in \partial _{\xi} \psi(t,0) } \langle q,y'(t)\rangle .
\end{displaymath}
  Let us fix $t \in [ a,b]$ such that these inequalities are
satisfied. The set $\partial _{\xi} \psi(t,0)$ being convex, we deduce that for
some $q \in
\partial _{\xi} \psi(t,0)$ we have $ \langle q,y'(t)\rangle= L (y
(t),y' (t)) - c$.
  For every $t \in [ a,b]$ let
\begin{equation}
   \label{dm368} B (t):= \left\{ q \in \partial_\xi \psi(t,0) \;|\;
\left\langle   q,y' (t)\right\rangle  =L (y (t),y' (t)) -c \right\}.
\end{equation}
By the previous argument  $B (t) \ne \emptyset $ for  a.e.\ $t \in
[ a,b]$. The graph of the set-valued map $t \leadsto B (t)$ is the
intersection of the sets  $B_1$ and $B_2$ defined by
\begin{eqnarray*}
& B_1:= \left\{  (t,q) \in [ a,b] \times \R^n\;|\; \left\langle   q,y'
(t)\right\rangle =L (y (t),y' (t))-c  \right\} ,
\\
  & B_2 :=
  \left\{  (t,q) \in [ a,b] \times \R^n\;|\;
  \forall \; \xi \in \R^n,\;\; \left\langle   q, \xi \right\rangle
\leq \psi (t, \xi )  \right\}.
\end{eqnarray*}
Clearly  $B_1$ and   $B_2 $ belong to $ {\cal L}_1 \times {\cal
B}_n$, where ${\cal B}_n$ denotes the $\sigma$-algebra
of all Borel subsets of $\R^n$. This implies that the graph of the
set-valued map $t
\leadsto B (t)$ defined by (\ref{dm368}) belongs to  ${\cal L}_1
\times {\cal B}_n$, and by the Aumann-von Neumann selection
theorem (see \cite[Theorem III.22]{cava}), there exists a
measurable function $p \colon  [ a,b] \to \R^n$ such that  $p (t)
\in B (t)$ for a.e.\ $t \in [ a,b]$. Then (\ref{dm362*}) and
(\ref{dm3a60*}) follow from (\ref{dm368}).  $\; \; \Box$

\vspace{5 mm}
Let $L_0:=\overline {\rm co}_u\, L$ be the lower semicontinuous
convex envelope
of  $L$ with respect to  $u$. Then $ (t,u) \mapsto L_0 (y (t),u)$ is a
Carath\'eodory function. This can be verified as in the case of $f_0$
(proof of Lemma  \ref{dm309}),
taking convex combinations of $n+1$ vectors.

     \begin{Theorem}
\label{dm3pr70}
Let $L\colon \R^n \times \R^n \rightarrow \R_+$ be a Borel function
and
let  $y \in W^{1,1} (a,b; \R^n)$ be a function which satisfies
(\ref{dm305}).
Suppose that for a.e.\ $t \in [ a,b]$
\begin{eqnarray}
   \label{ex1}
& L (y (t),y' (t))=L_0 (y (t),y' (t)),
\\
   \label{ex2}
   &
   -d^{ \,l}_v g_0 (t,1) \le d_u L_0 (y (t),y' (t)) (-y' (t)),
\\
   \label{ex3}
   &
    d^{ \,r}_v g_0 (t,1) \le d_u L_0 (y (t),y' (t)) (y' (t)) ,
\end{eqnarray}
where $d_u L_0(x,u)(\xi)$ denotes the
directional derivative of the convex function $L_0(x,\cdot)$ at $u$ in
the direction $\xi$.
Then there exist a constant $c \in \R$ and a measurable function $p
\colon [a,b] \rightarrow \R^n$ such that a.e.\ $t \in [ a,b]$
\begin{eqnarray}
   \label{dm3e550}
& p (t) \in \partial_u L_0 (y (t),y' (t)),
\\
   \label{dm3e560}
& L_0 (y (t),y' (t))- \left\langle   p (t),y' (t)\right\rangle =c.
\end{eqnarray}
Consequently,
\begin{equation}\label{constH0}
\langle p(t),y'(t)\rangle - L(y(t),y'(t)) =
\sup_{u\in\scriptR^n} (\langle p(t),u\rangle -L(y(t),u)) = -c
\end{equation}
for a.e.\ $t \in [ a,b]$.
\end{Theorem}

\begin{Remark}{\rm Since $g_0(t,v)\ge L_0(y(t),vy'(t))$ for a.e.\
$t\in [a,b]$ and for every $v\in\R$,
if (\ref{ex1}) holds, then by (\ref{dm3e46})
\begin{displaymath}
  -d^{ \,l}_v g_0 (t,1) \ge d_u L_0 (y (t),y' (t)) (-y' (t))
\quad\&\quad
    d^{ \,r}_v g_0 (t,1) \ge d_u L_0 (y (t),y' (t)) (y' (t)) ,
\end{displaymath}
so that  (\ref{ex2}) and  (\ref{ex3}) are actually equalities.
Assumptions (\ref{ex1})--(\ref{ex3}) are satisfied, for instance,
if $g_0(t,v)= L_0 (y (t),v y' (t))$
for a.e.\ $t \in [ a,b]$ and every  $v\in\R$.
}
\end{Remark}

\begin{Remark}{\rm If $H$ is the Hamiltonian associated to $L$,
defined in (\ref{H}), then,  by (\ref{constH0}),
$H(y(t),p(t))=-c$ for almost all $t\in [a,b]$.
The function $p$ corresponds to the co-state of optimal
control theory. In other words, we proved that the Hamiltonian
is constant along the optimal trajectory/co-state pair $(y,p)$. In the
case of smooth Hamiltonians this is
indeed a well known property of
optimal trajectories of autonomous
Bolza control problems.
}
\end{Remark}

\noindent
{\bf Proof of Theorem \ref{dm3pr70}} --- \hspace{ 2 mm}
The result follows from Lemma~\ref{dm363n}, taking
$\psi(t,\xi):=d_u L_0 (y (t),y' (t))(\xi)$. Indeed, the convexity
of $ L_0 (y (t),\cdot)$ implies that  $\psi(t,\cdot)$ is convex and
$\partial_\xi\psi(t,0)=\partial_u L_0 (y (t),y' (t))$
for a.e.\ $t \in [ a,b]$.
Equality (\ref{constH0}) follows from (\ref{ex1}), (\ref{dm3e550}),
and (\ref{dm3e560}).
   $\; \; \Box$

   \vspace{ 5mm}

To state further extensions of the DuBois-Reymond necessary condition, we
need to recall several notions of generalized derivatives.
Let $\varphi \colon \R^m\to \R \cup \{+ \infty \}$. The subdifferential
of $\varphi$ at $x\in \mbox{\rm dom}(\varphi)$ is defined by
\begin{equation}\label{subd3}
\partial^- \varphi(x):=\left\{
p\in \R^n\;|\; \liminf_{y\to x}\frac{ \varphi(y)- \varphi(x)-\langle
p,y-x\rangle}{|y-x|} \geq 0\right\} .
\end{equation}
An equivalent definition of subdifferential uses the lower
contingent derivatives of $\varphi$ defined by
\begin{equation}\label{dirder}
\forall \; u \in \R^m,\;\; D_{\uparrow } \varphi (x) (u) := \liminf_{
\textstyle {h \rightarrow 0+ \atop v \rightarrow u}} \frac{\varphi
(x+hv) - \varphi (x)}{h} .
\end{equation}
Then,
\begin{equation}
   \label{subd1}
   p \in \partial^- \varphi (x)  \;\; \Longleftrightarrow \;\; \forall \; v\in
\R^n, \;\; \left\langle   p,v\right\rangle \leq D_{\uparrow } \varphi
(x) (v) .
\end{equation}
(see, e.g., \cite{af90sva}). The upper contingent
derivative of $\varphi$ at $x$ is
defined by
\begin{displaymath}
\forall \; u \in \R^m,\;\; D _{\downarrow }  \varphi (x) (u) := \limsup_{
\textstyle {h \rightarrow 0+ \atop v \rightarrow u}} \frac{\varphi
(x+hv) - \varphi (x)}{h} .
\end{displaymath}
The superdifferential $\partial^+ \varphi (x)$ of $\varphi$ at $x$ is defined
by $ \partial^+ \varphi (x):=- \partial^- (-\varphi) (x)$ or, equivalently, by
\begin{equation}
   \label{supd1}
   p \in \partial^+ \varphi (x)  \;\; \Longleftrightarrow \;\; \forall \; v\in
\R^n, \;\; \left\langle   p,v\right\rangle \geq D_{\downarrow } \varphi
(x) (v) .
\end{equation}
We use also the lower Dini directional derivative, defined by
\begin{equation}
    \label{dm3e49}
\forall \; u \in \R^m,\;\; d^- \varphi(x) ( u ):=
\liminf_{h \rightarrow 0+} \frac{\varphi (x+h
u  ) -\varphi (x)}{h}.
\end{equation}

\vspace{ 5mm}

Let us return to the Lagrangian $L$ considered at the beginning of
this section. Partial derivatives and partial differentials of $L$
with respect to $u$
are defined in
the usual way: given $ x \in \R^n $,
we consider the function $\varphi (\cdot):=L(x,
\cdot)$, and set $ D _{\uparrow u} L (x,u) := D _{\uparrow } \varphi (u) $,
$
\partial _u^- L (x,u):=
\partial^- \varphi(u)$,
$ D _{\downarrow u} L (x,u) := D _{\downarrow } \varphi (u)$, $
\partial _u^+ L (x,u):= \partial^+ \varphi (u)$, and
$d^-_u L(x,u):= d^-\varphi(u)$.

\begin{Remark}\label{rem90}{\rm
By  (\ref{dm3e45}) and
(\ref{dm3e46}) for a.e.\ $t \in [ a,b]$ and every $v \in (0,2)$ we
have
\begin{displaymath}
   g_0 (t,v) \leq g (t,v) =L (y (t),vy' (t))\quad\&\quad g_0 (t,1)=g (t,1)=L (y
(t),y' (t)),
\end{displaymath}
which implies
\begin{equation}
        \label{dm364} \ \ \ \
-d^{\,l}_v g_0 (t,1) \leq  d_{u}^- L ( y (t),y' (t)) (-y' (t)) \quad\&
\quad
    d^{\,r}_vg_0 (t,1) \leq d_{u}^-L (y (t),y' (t)) (y' (t)).
\end{equation}
Therefore the conclusions of Lemma~\ref{dm363n} continue to hold if
(\ref{dm360'n}) is replaced by
\begin{equation}
        \label{dm360z}
        d_{u}^-L (y (t),y' (t)) (\pm y' (t))\le \psi(t, \pm y' (t))
\end{equation}
for a.e.\ $t \in [ a,b]$.
}\end{Remark}

\vspace{ 5mm}

For every $(x,u)\in
\R^n\times\R^n$ let
$\xi\!\mapsto \overline {\rm co} \, D _{\uparrow u} L(x,u)(\xi)$ be
the lower semicontinuous convex envelope of the function\
$\xi\mapsto D _{\uparrow u} L(x,u)(\xi)$.

\begin{Theorem}\label{dm363}
Let $L\colon \R^n \times \R^n \rightarrow \R_+$ be a Borel function and
let  $y \in W^{1,1} (a,b; \R^n)$ be a function which satisfies
(\ref{dm305}).
Suppose that for a.e.\ $t \in [ a,b]$
\begin{eqnarray}
\label{dm360'}
&
   -d^{ \,l}_v g_0 (t,1) \le \overline {\rm co} \,
D_{\uparrow u} L (y (t),y'(t))(-y' (t)),
\\
   \label{dm360'b}
   &
    d^{ \,r}_v g_0 (t,1) \le \overline {\rm co} \,
D_{\uparrow u} L (y (t),y'(t))(y' (t)),
\\
\label{dm360''}
& \forall \;\xi \in \R^n,\;\;\; \overline {\rm co}\, D _{\uparrow u}
L (y (t),y'(t))(\xi)\in\R.
\end{eqnarray}
Then there exist a constant  $c \in \R$ and a measurable function
$p \colon  [ a,b] \rightarrow \R^n$ such that
\begin{eqnarray}
   \label{dm362}
& p (t) \in \partial _{u} ^-L (y (t),y' (t)) \;\; \mbox{ for
a.e.\ } t \in [ a,b],
\\
\label{dm3a60}
&   L (y (t),y' (t))- \left\langle   p (t), y' (t)\right\rangle = c \;\;
\mbox{ for a.e.\ } t \in [ a,b].
\end{eqnarray}
\end{Theorem}

\begin{Remark}\label{rem105}{\rm By (\ref{dm364}) inequalities
(\ref{dm360'}) and (\ref{dm360'b}) are satisfied if
\begin{displaymath}
     d_{u}^-L (y (t),y' (t)) (\pm y' (t))\le \overline {\rm co} \,
D_{\uparrow u} L (y (t),y'(t))(\pm y' (t)),
\end{displaymath}
for a.e.\ $t \in [ a,b]$.
This shows that (\ref{dm360'})--(\ref{dm360''})
are always satisfied if $L(y(t),\cdot)$ is
differentiable at $y'(t)$.
}\end{Remark}

\begin{Remark}\label{rem106}{\rm
We recall that
a function $\varphi \colon  \R^n \to  \R$ is called {\em
semiconvex\/} if there exists $\omega \colon \R_{+} \times {\bf
R}_{+} \to  \R_{+}$, satisfying
          \begin{displaymath}\label{hp002}
          \forall \;  r \leq R,  \; \forall \;  s \leq S, \;  \; \omega (r,s)
\leq \omega (R,S) \;   \; \&   \; \;  \lim_{s\rightarrow 0+}\omega
(R,s) \;  =  \;  0,
          \end{displaymath}
such that for every $R>0$,  $\lambda \in [0,1]$, and all
          $ x, y \in B(0,R)$
          \begin{displaymath}
          \varphi (\lambda x + (1-\lambda )y) \leq \lambda \varphi (x)  +
(1-\lambda )\varphi (y)  +  \lambda (1-\lambda )|x-y| \; \omega
(R,|x-y|) .
          \end{displaymath}
Observe that every convex function  is semiconvex (with $\omega $
equal to zero). Furthermore, if   $\varphi \colon \R^{n}\to {\bf
R}$ is continuously differentiable, then it is semiconvex.
Using standard arguments of convex analysis (see, e.g.,
\cite[p.~25]{aubconvex})
one can prove that every semiconvex function $\varphi$ is
locally Lipschitz. Furthermore, for every $v  \in \R^n$
          \begin{displaymath}
          D _{\uparrow } \varphi (x) (v)  = \lim_{ h\rightarrow 0+}
         \frac{\varphi (x+hv ) -\varphi (x)}{h},
          \end{displaymath}
          and the function $v \mapsto  D _{\uparrow }
\varphi (x) (v)$ is convex (see, e.g., \cite[Theorem 3.9]{cf91cha}).
These facts, together with Remark~\ref{rem105}, show that
assumptions (\ref{dm360'})--(\ref{dm360''}) of Theorem \ref{dm363}
are always satisfied
when $L(x,\cdot)$ is semiconvex.
}\end{Remark}

  \noindent
{\bf Proof of Theorem \ref{dm363}} --- \hspace{ 2 mm} Let us define
$\psi(t,\xi):= \overline {\rm co} \,
D_{\uparrow u} L (y (t),y'(t))(\xi)$. Using the projection theorem
it is possible to check the measurability with respect to~$t$.
Notice that $\partial_\xi\psi(t,0)\subset \partial_u^- L (y (t),y'(t))$
for a.e.\ $t\in  [ a,b]$. Indeed, if
$q\in \partial_\xi\psi(t,0)$, then $\langle q,\xi\rangle\le
\psi(t,\xi)\le  D_{\uparrow u} L (y (t),y'(t))(\xi)$ for every
$\xi\in\R^n$, hence
$q\in \partial^-_u L (y (t),y'(t))$
by  (\ref{subd1}).
The conclusion follows then from
Lemma~\ref{dm363n}.
$\;\; \Box$

\begin{Remark}{\rm
Theorem~\ref{dm3pr70} has stronger assumptions and
stronger conclusions than Theorem~\ref{dm363}.
Indeed, as  $L \geq L_0$, it follows from
(\ref{ex1}) that
\begin{equation}\label{eq700}
     d_u L _0(y (t),y'(t))(\xi)=
D_{\uparrow u} L _0(y (t),y'(t))(\xi)\le
D_{\uparrow u} L (y (t),y'(t))(\xi)
\end{equation}
for every $\xi\in\R^n$. Since
$\xi\mapsto d_u L _0(y (t),y'(t))(\xi)$ is convex, we conclude that
\begin{displaymath}
d_u L _0(y (t),y'(t))(\pm y'(t))\le
  \overline {\rm co} \,
D_{\uparrow u} L (y (t),y'(t))(\pm y' (t)).
\end{displaymath}
This shows that
(\ref{ex1})--(\ref{ex3}) imply (\ref{dm360'})--(\ref{dm360''}).

On the other hand, (\ref{subd1}) and (\ref{eq700})
yield $ \partial ^- _u L_0 (y
(t),y' (t)) \subset \partial _u ^- L (y (t),y' (t))$ for a.e.\ $t
\in  [ a,b]$. Therefore  (\ref{ex1}), (\ref{dm3e550}), and
(\ref{dm3e560}) imply (\ref{dm362}) and (\ref{dm3a60}).}
\end{Remark}

\begin{Theorem}\label{dm363x}
Let $L\colon \R^n \times \R^n \rightarrow \R_+$ be a Borel function and
let  $y \in W^{1,1} (a,b; \R^n)$ be a function which satisfies
(\ref{dm305}).
Suppose that $L(x,\cdot)$ is locally Lipschitz continuous for every
$x\in\R^n$.
Then there exist a constant  $c \in \R$ and a measurable function
$p \colon  [ a,b] \rightarrow \R^n$ such that
for a.e.\ $t \in [ a,b]$
\begin{eqnarray}
   \label{dm362**}
& p (t) \in \partial_u L (y (t),y' (t)),
\\
\label{dm3a60**}
& L (y (t),y' (t))- \left\langle   p (t), y' (t)\right\rangle = c,
\end{eqnarray}
where $\partial_u L(x,u)$ denotes the Clarke generalized gradient of
$L(x,\cdot)$ at~$u$.
\end{Theorem}

\noindent
{\bf Proof} --- \hspace{ 2 mm}
Let us define
\begin{displaymath}
     \psi (t,\xi)= \limsup _{
\textstyle {h \rightarrow 0+ \atop u \rightarrow y'(t)}} \frac{L
(y (t), u+h\xi) - L (y (t),u)}{h} .
\end{displaymath}
It is known that $\psi (t,\cdot)$ is convex and that
  $\partial_\xi\psi (t,0)$ is the
Clarke generalized gradient of $L(y(t),\cdot)$ at $y'(t)$
  (see \cite{Cla}).
Since $d^-_{u} L (y (t), y' (t)) (\xi)\le \psi(t,\xi)$, the result
follows from
Lemma~\ref{dm363n} and Remark~\ref{rem90}.
$\;\;\Box$.

   \vspace{5 mm}

   Replacing subdifferential by superdifferential we get another
   extension of the DuBois-Reymond  necessary condition, which is
   meaningful only at those points $t \in [ a,b]$ for which
   $\partial _u^+ L (y (t), y' (t))\neq\emptyset$.

\begin{Proposition}
Let $L\colon \R^n \times \R^n \rightarrow \R_+$ be a Borel function
and
let  $y \in W^{1,1} (a,b; \R^n)$ be a function which satisfies
(\ref{dm305}).
There exists a constant  $c \in \R$ such that
\begin{displaymath}
  \forall\;   p \in \partial _u^+ L (y (t), y' (t)),\;\;\;
L (y (t),y' (t))- \left\langle   p,y' (t)\right\rangle =c
\end{displaymath}
for a.e.\ $t \in [ a,b]$.
\end{Proposition}

\noindent
{\bf Proof} --- \hspace{ 2 mm} {}From (\ref{dm3e46}) we have
\begin{eqnarray*}
& d^{\,r}_v g_0 (t,1) \leq D _{\downarrow u} L (y (t), y' (t)) (y' (t)),
\\
&  d^{\,l}_v g_0 (t,1) \geq -D _{\downarrow u} L (y (t), y' (t)) (-y'
   (t)).
\end{eqnarray*}
These inequalities and
(\ref{dm3e47}) imply that
there exists a constant  $c$ such that for a.e.\ $t \in [
a,b]$
\begin{displaymath}
-D _{\downarrow u} L (y (t), y' (t)) (-y'
   (t)) \leq L (y (t),y' (t)) -c \leq
D _{\downarrow u} L (y (t), y' (t)) (y' (t)),
\end{displaymath}
and
we deduce from (\ref{supd1}) that for all $p \in \partial _u^+ L
(y (t), y' (t))$
\begin{displaymath}
   \left\langle -p, -y' (t)  \right\rangle \leq L (y (t),y' (t)) -c  \leq
\left\langle   p, y' (t)\right\rangle,
\end{displaymath} ending the proof.
$\; \; \Box$

\section{Hamilton-Jacobi Inequalities}\label{Hamilton}

Let  $ \varphi \colon \R^n \to \R_+ \cup \{ + \infty \}$ with
$\varphi \not\equiv+ \infty $.
Given $T>0$ and $y_0 \in \R^n$, let us consider the
Bolza problem:
\begin{displaymath}
\min_y \int_{0}^{T} L (y(s), y' (s))ds + \varphi  (y(T))
\end{displaymath}
over all absolutely continuous functions $y \in W^{1,1} (0,T;
\R^n)$ satisfying the initial condition $y (0)= y_0$. The dynamic
programming approach associates with this problem the family of
problems ($t\geq 0$, $x\in\R^n$):
\begin{displaymath}
\min_y \int_{0}^{t} L (y (s), y' (s))ds + \varphi  (y(t))
\end{displaymath}
over all absolutely continuous functions $y \in W^{1,1} (0,t;
\R^n)$ satisfying $y (0)= x$. The corresponding value function
$V\colon  \R_+
\times \R^n \to  \R_+ \cup \left\{ + \infty  \right\}$ is defined
by (\ref{value}).

\begin{Proposition}\label{prop15}
Let $L\colon \R^n \times \R^n \rightarrow \R_+$ be a Borel function
and let $ \varphi \colon \R^n \to \R_+ \cup \{ + \infty \}$ with
$\varphi \not\equiv+\infty$.
Then $V(0,x)=\varphi(x)$ for every $x\in\R^n$.
Furthermore, if $L$ satisfies   (\ref{dm302}) and (\ref{dm304}), and
$\varphi$ is lower semicontinuous, then
\begin{equation}\label{eq115}
\liminf_{
\textstyle {t \rightarrow 0+ \atop x \rightarrow x_0}} V(t,x)\ge
\varphi(x_0)
\end{equation}
for every $x_0\in\R^n$.
\end{Proposition}

\noindent
{\bf Proof} --- \hspace{ 2 mm} Fix $x _0\in \R^n$ and
let $t_i \to 0+$, $x_i \to x _0$ be such that
\begin{displaymath}
\liminf_{
\textstyle {t \rightarrow 0+ \atop x \rightarrow x_0}} V(t,x)
= \lim_{i \to \infty} V(t_i,x_i).
\end{displaymath}
If the above limit is infinite, then our claim follows. If this limit
is finite,
then we consider $y_i \in W^{1,1}(0,t_i;\R^n)$ such
that $y_i(0)=x_i$ and
\begin{displaymath}
\int_{0}^{t_i} L (y_i (s), y'_i(s))ds +
\varphi  (y_i(t_i)) \leq V(t_i,x_i) + \frac{1}{i}.
\end{displaymath}
By (\ref{dm304}), since $\varphi \geq 0$,
for some $M>0$ we have $\int_{0}^{t_i} \Theta( y'_i(s))ds \leq M$ for
every~$i$. Setting $y_i'(s)=0$ for
  $s \in (t_i,1]$, we deduce from (\ref{dm302}) that the functions $y_i'$ are
  equiintegrable and therefore the functions $y_i$ are equicontinuous.
  Since $t_i \to 0+$ and $y_i(0)=x_i\to x_0$, we get $y_i(t_i)\to x_0$.
On the other hand, since $L \geq 0$,  we have
$\varphi  (y_i(t_i)) \leq V(t_i,x_i) + \frac{1}{i}$.
Taking the lower limit and using
the lower semicontinuity of $\varphi$ we conclude the proof.
$\;\;\Box$

\vspace{5 mm}
We recall that $\R_+ ^{\star}:=\{t\in\R \;|\; t>0\}$.
In this section we often assume the following hypotheses:

\vspace{ 5 mm}

\noindent
({\bf H1}) for every $(t,x) \in \R_+^\star \times \R^n$ the infimum in
(\ref{value}) is attained,

\medskip

\noindent
({\bf H2}) $L$ is  locally bounded and satisfies
(\ref{dm302}) and~(\ref{dm304}).

\vspace{ 5 mm}

\noindent
It is easy to see that ({\bf H2}) implies that
$0 \leq V(t,x) < \infty$ for all $(t,x) \in \R_+^\star \times \R^n$.

\begin{Remark}\label{rem10}
{\rm If for every $t>0$ the functional
\begin{equation}\label{eq50}
{\cal L}_0^t(y) :=
\int_{0}^{t} L(y(s),y'(s))ds
\end{equation}
is  sequentially weakly lower
semicontinuous  on $W^{1,1}(0,t;\R^n)$ and $ \varphi$ is lower
semicontinuous, then from (\ref{dm302}) and (\ref{dm304})
it follows that  ({\bf H1})
is satisfied. Furthermore, arguing as in \cite[Proof of
Proposition 3.1]{dm-fra}, we can show that in this case $V$ is lower
semicontinuous on $\R_+ \times \R^n$.
}
\end{Remark}

\begin{Lemma}\label{minlip}
Let $L\colon \R^n \times \R^n \rightarrow \R_+$ be a Borel function
and let $ \varphi \colon \R^n \to \R_+ \cup \{ + \infty \}$ with
$\varphi \not\equiv+\infty$. Assume that $({\bf H1})$ and $({\bf
H2})$ are satisfied.
Then, given $ (t_0,x_0) \in \R_+^{ \star } \times \R^n$
and $0<\delta<t_0$,
there exists $r>0$ such that for all $ (t,x) \in B ( (t_0,x_0),
\delta )$  every
minimizer $ y ( \cdot; t, x)$ of  (\ref{value})  is $r$-Lipschitz.
\end{Lemma}

\noindent
{\bf Proof} --- \hspace{ 2 mm} Consider $y_0 \in \R^n$, with
$\varphi(y_0)<+\infty$,
and set $z(s)= x+ \frac{s}{t}(y_0-x)$. If $y(\cdot; t,x)$ is a
minimizer of (\ref{value}), we obtain
\begin{displaymath}
\int _{0}^{t}
L(y(s;t,x),y'(s;t,x))ds \leq \varphi (y_0) + t \sup_{s \in [0,t]}
L (z(s),(y_0-x)/t ).
\end{displaymath}
Since $L$ is locally bounded, for every
$0<\delta<t_0$
there exists a constant $M_\delta>0$ such that
\begin{displaymath}
\sup_{s \in [0,t]} L (z(s),(y_0-x)/t )\le M_\delta
\end{displaymath}
for every   $ (t,x) \in B ( (t_0,x_0), \delta )$. The conclusion follows
now from Theorem~\ref{dm306a}.
   $\; \; \Box$

\begin{Theorem}\label{liploc}
Let $L\colon \R^n \times \R^n \rightarrow \R_+$ be a Borel function
and let $ \varphi \colon \R^n \to \R_+ \cup \{ + \infty \}$ with
$\varphi \not\equiv+\infty$. Assume that $({\bf H1})$ and $({\bf
H2})$ are satisfied.  Then $V $ is locally Lipschitz on
$\R_+^{ \star } \times \R^n$. If $\varphi$ is lower semicontinuous
on $\R^n$, then $V $ is lower semicontinuous  on
${\R_+\times \R^n}$.
\end{Theorem}

\noindent
{\bf Proof} --- \hspace{ 2 mm}
The Lipschitz
continuity is proved in \cite[Corollary 3.4]{dm-fra}.
For the reader's convenience we repeat here the proof.

  Fix
  $(t_0,x_0) \in \R_+^{ \star } \times \R^n$.
By Lemma
\ref{minlip}, there exist $r>0$ and $\delta >0$ such that for all $ (t,x)
\in
B ( (t_0,x_0), \delta )$  every minimizer $ y ( \cdot; t, x)$ of
(\ref{value}) is $r$-Lipschitz. We may assume that $5\delta<t_0$.
Let $(t_1, x_1)$ and $(t_2, x_2)$ be two distinct points of
$B ( (t_0,x_0), \delta )$, let $h_1:=|t_1-t_2| + |x_1-x_2|$, and
$s_1:= h_1-t_1+t_2$. Let $u_1\in\R^n$ be
such that $y(s_1;t_2,x_2)= x_1+h_1u_1$. Then $0< h_1<t_1$,
$0\le s_1 \le 2h_1$, and
\begin{equation}\label{dm301*}
|u_1| \leq \frac{|y(s_1;t_2,x_2)-x_2|}{h_1} +
\frac{|x_2-x_1|}{h_1} \leq 2r + 1.
\end{equation}
Let $y_1\colon [0,t_1]\mapsto\R^n$ be the function defined by
\begin{displaymath}
  y_1 (s)=\cases{x_1+su_1&if $0\le s\le h_1$,
  \cr
  \noalign{\vskip 3pt}
  y(s-t_1+t_2; t_2,x_2)&if $h_1\le s \le t_1$.
  \cr}
\end{displaymath}
Then
\begin{eqnarray*}
& \displaystyle
V(t_1,x_1)\leq \int_0^{t_1} L(y_1(s), y_1^\prime(s))ds +
\varphi(y_1(t_1))
=
\\
& \displaystyle
=\int_0^{h_1} L(x_1+s u_1, u_1)ds
+
\int_{s_{1}}^{t_{2}} L(y(s; t_2,x_2), y'(s; t_2,x_2))ds
+\varphi(y(t_{2}; t_2,x_2)).
\end{eqnarray*}
As $s_{1}=h_{1}-t_{1}+t_{2}\geq 0$ and $L\geq 0$, we obtain
\begin{displaymath}
V(t_1,x_1) \leq \int_0^{h_1} L(x_1+s u_1, u_1)ds +
V(t_2,x_2).
\end{displaymath}
Since $L$ is locally bounded, it follows from (\ref{dm301*}) that
there exists a constant $M$, depending only on $L$, $t_0$, $x_0$,
$\delta$, and $r$, such that
\begin{displaymath}
V(t_1,x_1)-V(t_2,x_2) \leq Mh_1=M(|t_1-t_2| + |x_1-x_2|).
\end{displaymath}
Exchanging the roles of $(t_1,x_1)$ and $(t_2,x_2)$ we get that
$V$ in $M$-Lipschitz on $B ((t_0,x_0), \delta )$.

If $\varphi$ is lower semicontinuous,  then $V$
is lower semicontinuous at all points of $\{0\}\times\R^n$ by
Proposition~\ref{prop15}. The lower semicontinuity on
$\R_+^\star\times\R^n$ is a consequence of the local
Lipschitz continuity.
$\; \; \Box$

\vspace{ 5 mm}

We recall the definition of the function $L ^+(x,u)$ used in
\cite{dm-fra} to study Hamilton-Jacobi equation (\ref{HJ}) arising
from a
discontinuous Lagrangian:
\begin{equation}\label{feb30}
L^+(x,u) := \limsup_{h \rightarrow 0+} \frac{1}{h} \inf_y \left\{
\int_{-h}^{0} L (y (s),y' (s))ds  \;|  \; y (-h)=x-hu,\; y (0)=x\right\}.
   \end{equation}
Similarly we define the function $L^-(x,u)$ by
\begin{equation}\label{dm384}
L^-(x,u) := \liminf_{h \rightarrow 0+} \frac{1}{h} \inf_y \left\{
\int_{0}^{h} L (y (s),y' (s))ds   \;|  \; y (h)=x+hu,\; y (0)=x\right\} .
   \end{equation}

\begin{Proposition} \label{eq---}
Assume that $L$ is locally bounded. If $u_h\to u$ as $h\to
0+$, then
\begin{eqnarray}\label{eq20}
&\displaystyle
L^+ (x,u) =\limsup_{h \rightarrow 0+}
\frac{1}{h} \inf _y\left\{ \int_{-h}^{0} L (y ( s ),y' ( s )) d s
\;|\; y (-h)=x-hu_h  ,\;  y (0) = x\right\},\\
\label{eq21}
&\displaystyle
L^- (x,u) =\liminf_{h \rightarrow 0+}
\frac{1}{h} \inf _y\left\{ \int_{0}^{h} L (y ( s ),y' ( s )) d s
\;|\; y (h)=x+hu_h  ,\;  y (0) = x\right\}
.
\end{eqnarray}
\end{Proposition}

\noindent
{\bf Proof} --- \hspace{ 2 mm}
The following proof is a slight modification of the proof of
\cite[Proposition~3.6]{dm-fra}.
Let us fix $ (x,u)$ and $u_h$ as in the statement of the
proposition, and let $ \overline L (x,u)$ be the right hand side
of  (\ref{eq20}).
We want to show that $\overline L(x,u)\leq L^+(x,u)$.
For every $h>0$ let $\varepsilon_{h}= |u_h-u|$
and let $y_{h}$ be such that $y_h(-(1-\varepsilon_{h})h)
=x-(1-\varepsilon_{h})h u $,
  $y_h (0)=x$, and
\begin{eqnarray*}
   &\displaystyle
   \int_{-(1-\varepsilon_{h})h}^{0} L(y_h ( s ),y'_h ( s )) d s
   -h^2\le
   \\
     &\displaystyle
  \le \inf_y \left\{ \int_{-(1-\varepsilon_{h})h}^{0}
  L(y ( s ),y' ( s )) d s  :
y(-(1-\varepsilon_{h})h)=x-(1-\varepsilon_{h})h u  ,\;
  y (0)=x\right\}.
\end{eqnarray*}
We extend  $y_h$ on the interval
$ [-h,-(1-\varepsilon_{h})h]$ by the
affine function satisfying $y_h (-h)= x-h u_h$ and
$y_h(-(1-\varepsilon_{h})h)=x-(1-\varepsilon_{h})h u$.
Since on this interval the derivative of $y_h$ is equal to
$ (u_h- (1-\varepsilon_{h}) u)/\varepsilon_{h}$, which is uniformly
bounded,
we deduce that for some $M>0$ and all $h>0$,
\begin{displaymath}
  \int_{-h}^{ 0} L (y_h (s),y_h' (s))ds \leq
\int_{-(1-\varepsilon_{h}) h}^{0}L (y_h
(s),y_h' (s))ds + \int_{-h}^{ -(1-\varepsilon_{h}) h} M ds .
\end{displaymath}
Dividing by $h$ and taking the upper limit
when $h \rightarrow 0+$
we get $\overline L (x,u) \leq L^+ (x,u)$. The opposite inequality can
be proved in the same way.
The proof of (\ref{eq21}) is similar.
$\; \; \Box$

  \begin{Remark} \label{rem1}{\rm
{}From the previous proposition it follows that, if $h_i\to 0+$,
$u_i\to u$, and
$y_i\in W^{1,1}(0,h_i)$ satisfies $y_i(0)=x$ and $y(h_i)=x+h_iu_i$, then
\begin{equation}\label{eq41}
L^-(x,u) \le \liminf_{i \rightarrow \infty}
\frac{1}{h_i} \int_{0}^{h_i} L (y_i ( s ),y'_i ( s )) d s.
\end{equation}
}
\end{Remark}

\begin{Proposition} \label{prop2}
Assume that $L$ is locally bounded.
For every $y\in W^{1,1}(0,t;\R^n)$
we have $L^+(y(s),y'(s)) \le L(y(s),y'(s))$
and $L^-(y(s),y'(s)) \le L(y(s),y'(s))$
for a.e.\ $s\in [0,t]$. If $y$ is a minimizer of (\ref{value}), then
$L^+(y(s),y'(s))=L^-(y(s),y'(s))=L(y(s),y'(s))$ for a.e.\ $s\in [0,t]$.
\end{Proposition}

\noindent
{\bf Proof} --- \hspace{ 2 mm}
Assume first $y\in W^{1,\infty}(0,t;\R^n)$.
Since $L$ is locally bounded, the function $s \mapsto
\psi (s):= \int_{0}^{s} L (y (\tau),y' (\tau))d\tau$ is absolutely
continuous. Let
$s\in [0,t]$ be such that both $ \psi ' (s)$ and $y' (s)$ do exist and
$\psi ' (s) = L ( y (s),y' (s)) $. Since
$u_{h}=(y(s)-y(s-h))/h$ converges to $y'(s)$ as $h\to 0+$,
from Proposition~\ref{eq---} we obtain
\begin{displaymath}
L^+(y(s),y' (s)) \leq \lim_{h\to 0+} \frac{1}{h} \int_{s-h}^{s} L ( y
(\tau),y' (\tau))
d\tau = \psi ' (s) = L ( y (s),y'(s)) ,
\end{displaymath}
which concludes the proof of the inequality
$L^+(y(s),y'(s)) \le L(y(s),y'(s))$ when $y$ is Lipschitz.

If $y\in W^{1,1}(0,t;\R^n)$, we can apply a Lusin
type approximation theorem for Sobolev functions (see, e.g.,
\cite[Theorem~3.10.5]{Zie}), which asserts that for every
$\varepsilon>0$ there exist $y_\varepsilon\in W^{1,\infty}(0,t;\R^n)$
and an open set $U_\varepsilon$ such that
$|U_\varepsilon|<\varepsilon$ and $y_\varepsilon(s)=y(s)$ for all
$s\in [0,t]\setminus U_\varepsilon$. As $y'_\varepsilon(s)=y'(s)$ for
a.e.\  $s\in [0,t]\setminus U_\varepsilon$, we obtain that
$L^+(y(s),y'(s)) \le L(y(s),y'(s))$
and $L^-(y(s),y'(s)) \le L(y(s),y'(s))$
for a.e.\ $s\in [0,t]\setminus U_\varepsilon$. Since $\varepsilon$ is
arbitrary, these inequalities hold for a.e.\ $s\in [0,t]$.

If $y$ is a minimizer of (\ref{value}), then for every $s\in (0,t)$
and every $h\in (0,s)$ we have
\begin{equation}\label{eq400}
\ \ \ \ \ \ \ \ \ \ \int_{s-h}^{s} L (y ( \tau ),y' ( \tau )) d \tau =
\inf _z\left\{ \int_{-h}^{0} L (z (\tau),z' (\tau)) d\tau
\;|\; z (-h)\!=\!y(s-h)  ,\;  z (0)\!= \!y(s)\right\}.
\end{equation}
Let us fix a Lebesgue point $s\in (0,t)$
for the function $L(y(s),y'(s))$
such that $(y(s)-y(s-h))/h\to y'(s)$ as $h\to 0+$. If we divide  both sides of
(\ref{eq400}) by~$h$, the left hand side tends to $L(y(s),y'(s))$
while the
right hand side tends to $L^+(y(s),y'(s))$ thanks to
Proposition~\ref{eq---}, applied with $x:=y(s)$, $u:=y'(s)$, and
$u_h:=(y(s)-y(s-h))/h$. Therefore
$L^+(y(s),y'(s))=L(y(s),y'(s))$ for a.e.\ $s\in [0,t]$.
The proof for $L^-(y(s),y'(s))$ is
similar.
$\; \; \Box$

\vspace{5 mm}

   Let us define
\begin{eqnarray}
\label{eq200}
&\displaystyle H ^+(x,p):= \sup_{u \in \scriptR^n} \left(
\left\langle   p,u\right\rangle
   -L^+(x,u) \right),
\\
\label{eq201}
&
\displaystyle
   H ^-(x,p):= \sup_{u \in \scriptR^n} \left( \left\langle   p,u\right\rangle
   -L^-(x,u) \right) .
\end{eqnarray}

\begin{Theorem}\label{thm2}
     Let $L\colon \R^n \times \R^n \rightarrow \R_+$ be a Borel function
and let $ \varphi \colon \R^n \to \R_+ \cup \{ + \infty \}$ with
$\varphi \not\equiv+\infty$. Assume that $({\bf H1})$ and $({\bf
H2})$ are satisfied.
Then the value function $V$ satisfies the following two
inequalities:
\begin{eqnarray} \label{dm360}
&\forall \; (t,x) \in \R_+^{ \star } \times \R^n ,\;\; \exists \; u
   \in \R^n,\;\;  D _{\uparrow}V (t,x) (-1,u)  \leq -L^- (x,u),
\\
\label{dm360c}
   &\forall \; (t,x) \in \R_+ \times \R^n , \;
   \forall \; u \in \R^n,\;\; D _{\downarrow}V (t,x) (1,-u) \leq L^+
   (x,u).
\end{eqnarray}
   Consequently, $V$ is a supersolution to (\ref{HJ}) on
$\R_+^{ \star } \times \R^n $  when $H$ is replaced by $H^-$,
i.e.,
\begin{equation} \label{supersol}
   \forall \; (t,x) \in \R_+^{ \star } \times \R^n ,\;\; \forall \; (p_t,p_x)
\in
\partial^-V (t,x),\;\;p_t+ H^-(x,-p_x) \geq 0,
\end{equation}
   and $V$ is a subsolution to (\ref{HJ}) on
$\R_+ \times \R^n $  when $H$ is replaced by $H^+$, i.e.,
\begin{equation} \label{dm360d}
   \forall \; (t,x) \in \R_+ \times \R^n , \;
   \forall \; (p_t, p_x) \in \partial^- V (t,x),\;\; p_t + H^+ (x,-p_x)  \leq
   0 .
\end{equation}
\end{Theorem}

\noindent
{\bf Proof} --- \hspace{ 2 mm}
   Let  $t>0, \; x \in \R^n$, and let $y$ be a
minimizer of  (\ref{value}).  By Theorem \ref{dm306} $y ( \cdot
)$ is Lipschitz. By minimality for all $0<h \leq t$ we have
\begin{displaymath}
   V (t,x)= V (t-h, y (h))+ \int_{0}^{h}L (y (s), y' (s))ds .
\end{displaymath}
Consider $h_i \rightarrow 0+$ such  that for some $u \in \R^n$,
$u_i:=(y (h_i)-x)/h_i\rightarrow u$. Then $y(h_i)=x+h_iu_i$, and
(\ref{eq41}) yields
\begin{displaymath}
   D _{\uparrow }V (t,x) (-1,u) \leq -\limsup_{i \rightarrow  \infty }
\frac{1}{h_i} \int_{0}^{h_i} L (y (s), y' (s))ds \leq -L^-(x,u) ,
\end{displaymath}
which proves (\ref{dm360}).

Let  $ (p_t,p_x) \in \partial^- V (t,x)$. Then, by (\ref
   {subd1}),
\begin{displaymath}
- p_t + \left\langle   p_x, u\right\rangle \leq
    D _{\uparrow }V (t,x) (-1,u),
\end{displaymath}
hence  $p_t + \left\langle   -p_x, u\right\rangle -L^- (x,u) \geq
0$. By (\ref{eq201}) this inequality gives (\ref{supersol}).

To prove inequality (\ref{dm360c}), we fix any $u \in \R^n$ and
let $h_i \rightarrow 0+$, $u_i \rightarrow u$.
{}From the definition of $V$ it follows that
\begin{eqnarray*}
&V(t+h_i,x-h_iu_i)-V (t,x) \leq
\\
&\displaystyle\le  \inf_y \left\{ \int _{-h_i}
^{0}L(y(s),y'(s))ds \;|\;  y(-h_i)=x-h_iu_i, y(0)=x \right\}.
\end{eqnarray*}
Then we divide by $h_i$ and pass
to the upper limit as $i\to\infty$. Taking (\ref{eq20})
 into account we obtain
(\ref{dm360c}). To prove (\ref{dm360d}) it is enough to apply
(\ref{subd1}), (\ref{dm360c}), and~(\ref{eq200}). $\; \; \Box$

\begin{Theorem}\label{thm1}
Let $L\colon \R^n \times \R^n \rightarrow \R_+$ be a Borel function,
let $ \varphi \colon \R^n \to \R_+ \cup \{ + \infty \}$
be a lower semicontinuous
function with
$\varphi \not\equiv+\infty$,  and
let $W\colon  \R_+
\times \R^n \to \R \cup \{ + \infty \}$ be a lower semicontinuous
function which satisfies the initial condition $W (0, \cdot )=
\varphi$. Assume that $({\bf H1})$ and $({\bf
H2})$ are satisfied. If $W$ is a subsolution of the Hamilton-Jacobi equation
(\ref{HJ}), in the sense that
\begin{equation} \label{feb45}
   \forall \; (t,x) \in \mbox{\rm dom}(W), \;
   \forall \; u \in \R^n,\;\; D _{\downarrow }W (t,x) (1,-u) \leq L^+ (x,u) ,
\end{equation}
then $W\le V$ on $\R_+ \times \R^n$.
Therefore the value function $V$ is the greatest
lower semicontinuous function $W$ which satisfies
inequality~(\ref{feb45}) and the initial condition
$W(0,\cdot)=\varphi$.
\end{Theorem}

We shall use the following well known lemma (see, e.g.,
\cite[Chapter~5, Section~2, Exercise~3]{Roy}). For completeness we
give here an elementary proof. To simplify notation, in the case of
functions $f$ of one real variable the lower Dini derivative $d^-f(t)(1)$,
defined by (\ref{dm3e49}), is denoted by $d^-f(t)$.

\begin{Lemma}\label{lemma1}
Let $[a,b]$ be a bounded closed interval in $\R$ and let
$f\colon [a,b]\to\R\cup\{+\infty\}$ be a lower semicontinuous
function such that $d^-f(t)\le 0$ for every $t\in [a,b)$ with
$f(t)<+\infty$.
Then $f(b)\le f(a)$.
\end{Lemma}

\noindent {\bf Proof} --- \hspace{ 2 mm} For every $\varepsilon>0$
let us consider the lower semicontinuous function
$f_\varepsilon(t):= f(t)-\varepsilon t$, which  satisfies the
stronger inequality $d^-f_\varepsilon(t)\le -\varepsilon$ for
every $t\in [a,b)$ with $f_\varepsilon(t)<+\infty$. We claim that
$f_\varepsilon(b)\le f_\varepsilon(a)$. If this inequality is not
satisfied for some $\varepsilon>0$, then the infimum of
$f_\varepsilon$ is attained at some $t_\varepsilon\in [a,b)$ and
we have $f_\varepsilon(t_\varepsilon)<+\infty$. Since
$d^-f_\varepsilon(t_\varepsilon)\le -\varepsilon$, there exists
$s_\varepsilon\in (t_\varepsilon, b)$ such that
$f_\varepsilon(s_\varepsilon)< f_\varepsilon(t_\varepsilon)$,
contradicting the minimality of $t_\varepsilon$ and proving our
claim. Taking the limit in the inequality $f_\varepsilon(b)\le
f_\varepsilon(a)$ as $\varepsilon\to0+$ we conclude the proof. $\;
\; \Box$

\vspace{5 mm}

An alternative proof of the following corollary can be found in
\cite[Chapter~6, Exercise~10]{Let}.

\begin{Corollary}\label{cor1}
Let $[a,b]$ be a bounded closed interval in $\R$ and let
$f\colon [a,b]\to\R\cup\{+\infty\}$ be a lower semicontinuous
function with $f(a)<+\infty$.
Suppose that there exists a constant $M\in \R$ such that
$d^-f(t)\le M$
for every $t\in [a,b)$ with $f(t)<+\infty$. Then
$f(t)<+\infty$ for every $t\in [a,b]$.
Suppose, in addition, that for some $g\in L^1(a,b)$
we have
$d^-f(t)\le  g(t)$
for a.e.\ $t\in [a,b]$.
Then
\begin{equation}\label{eq5}
f(b)-f(a) \le
\int_a^b g(t)dt.
\end{equation}
\end{Corollary}

\noindent
{\bf Proof} --- \hspace{ 2 mm} By Lemma~\ref{lemma1} the function
$f_M(t):=f(t)-Mt$ is nonincreasing. Therefore $f(t)\le
f(a)+M(t-a)<+\infty$ for every
$t\in[a,b]$. If $d^-f(t)\le  g(t)$
for a.e.\ $t\in [a,b]$, then $f'_M(t)\le g(t)-M$ for
a.e.\ $t\in [a,b]$, and the Lebesgue theorem on derivatives of
monotone functions yields
\begin{displaymath}
f_M(b)-f_M(a)\le \int_a^b f'_M(t)dt \le \int_a^b (g(t)-M)dt,
\end{displaymath}
which implies~(\ref{eq5}).
$\; \; \Box$

\vspace{5 mm}

\noindent
{\bf Proof of Theorem~\ref{thm1}} --- \hspace{ 2 mm} Let us fix
$t>0$, $x\in \R^n$, and let $y$ be
a minimizer of (\ref{value}).
It is Lipschitz continuous by Theorem~\ref{dm306}. Let us define
$\gamma(s):=W(s,y(t-s))$. Then $\gamma$ is lower semicontinuous on
$[0,t]$ and $\gamma(0)=\varphi(y(t))<+\infty$.
Let us fix $s\in [0,t)$ with $\gamma(s)<+\infty$.
Consider a sequence $h_i\to0+$
such that
\begin{equation}\label{eq6}
d^-\gamma(s)=
\lim_{i\to \infty} \frac{\gamma(s+h_i)-\gamma(s)}{h_i}.
\end{equation}
We can write
\begin{equation}\label{eq7}
\gamma(s+h_i)-\gamma(s)= W(s+h_i, y(t-s)-h_i u_i)-W(s, y(t-s)),
\end{equation}
with $u_i:=-(y(t-s-h_i)-y(t-s))/h_i$. Passing to a subsequence we may
assume that $u_i$ converges in $\R^n$ to some vector $u$, whose norm
is bounded by the
Lipschitz constant of $y$. {}From (\ref{feb45}), (\ref{eq6}), and (\ref{eq7})
it follows that
\begin{equation}\label{eq8}
d^-\gamma(s)\le
L^+(y(t-s),u).
\end{equation}
Since the function $L^+$ is locally bounded, we conclude that there
exists a constant $M$ such that
$
d^-\gamma(s)\le M
$
for every $s\in [0,t)$ with $\gamma(s)<+\infty$.
By Corollary~\ref{cor1} this implies that $\gamma(s)<+\infty$
for every $s\in [0,t]$.

If the derivative $y'(t-s)$ exists, then $u=y'(t-s)$. Therefore (\ref{eq8})
gives
$
d^-\gamma(s)\le
L^+(y(t-s),y'(t-s))$ for a.e.\ $s\in[0,t]$,
which, together with Proposition~\ref{prop2}, yields
$
d^-\gamma(s)\le
L(y(t-s),y'(t-s))
$ for a.e.\ $s\in[0,t]$.
By Corollary~\ref{cor1}
we obtain
\begin{displaymath}
\gamma(t)\le \gamma(0) + \int_0^t
L(y(t-s),y'(t-s))ds,
\end{displaymath}
which is equivalent to
\begin{displaymath}
W(t,x)\le
\varphi(y(t)) + \int_0^t L(y(s),y'(s))ds.
\end{displaymath}
Since the right
hand side is equal to $V(t,x)$, we have proved that $W(t,x)\le V(t,x)$.

The last assertion of the theorem follows now
from Theorems~\ref{liploc} and~\ref{thm2}.
$\; \; \Box$

\vspace{5 mm}

We conclude this section with some results which connect the
minimizers of (\ref{value}) with the contingent derivatives of the
value function.

\begin{Theorem}\label{dm391}
     Let $L\colon \R^n \times \R^n
\rightarrow \R_+$ be a Borel function and
let $ \varphi \colon \R^n \to \R_+ \cup \{ + \infty \}$ with
$\varphi \not\equiv+\infty$. Assume that $({\bf H1})$ and $({\bf
H2})$ are satisfied.
If $y$ is a minimizer of (\ref{value}),
then
\begin{eqnarray} \label{dm360x}
& D _{\uparrow}V (t-s,y(s)) (-1,y'(s)) =
D _{\downarrow}V (t-s,y(s)) (-1,y'(s))  = -L(y(s),y'(s)),
\\
\label{dm360y}
   & D _{\uparrow}V (t-s,y(s)) (1,-y'(s))=
   D _{\downarrow}V (t-s,y(s)) (1,-y'(s)) = L(y(s),y'(s))
\end{eqnarray}
for almost all $s \in [0,t]$.
\end{Theorem}

\noindent
  {\bf Proof} --- \hspace{ 2 mm}
Since $V$ is locally
Lipschitz on $\R_+^{ \star } \times \R^n$ by Theorem~\ref{liploc},
the function $\gamma(s):=V(t-s, y (s))$ is locally absolutely
continuous on $[0,t)$.
Fix $s \in (0,t)$ such
that the derivatives $\gamma'(s)$ and $y'(s)$ exist.
Let us prove that
\begin{eqnarray}
\label{eq502}
& D _{\uparrow}V (t-s,y(s)) (-1,y'(s))=
D _{\downarrow}V (t-s,y(s)) (-1,y'(s))=
\gamma'(s),
\\
\label{eq503}
& D _{\uparrow}V (t-s,y(s)) (1,-y'(s))=
  D _{\downarrow}V (t-s,y(s)) (1,-y'(s))=-
\gamma'(s).
\end{eqnarray}
Let $h_i\to0+$ and $u_i\to y'(s)$ such that
\begin{equation}\label{eq520}
\ \ \ \ \  D _{\downarrow}V (t-s,y(s)) (-1,y'(s))=
\lim_{i\to\infty} \frac{V(t-s-h_i, y(s)+h_iu_i)-
V(t-s, y(s))}{h_i} .
\end{equation}
As $V$ is Lipschitz near $(t-s, y(s))$, there exists a constant
$M>0$ such that
\begin{eqnarray}\label{eq521}
& |V(t-s-h_i, y(s+h_i))- V(t-s-h_i, y(s)+h_iu_i)|\le
\\
\nonumber
&\le M|y(s+h_i)-y(s)-h_iy'(s)|  + Mh_i|u_i-y'(s)|
\end{eqnarray}
for $i$ large enough. {}From (\ref{eq520}) and (\ref{eq521})
we obtain the second equality in (\ref{eq502}). The other equalities
in (\ref{eq502}) and (\ref{eq503}) can be obtained in the same way.

On the other hand, by minimality, for all small $h>0$
we have
\begin{equation}
\label{eq500} \displaystyle V(t-s-h, y (s+h))=V(t-s,y(s)) -
\int_{s}^{s+h}L(y(\tau),y'(\tau))d\tau .
\end{equation}
Since $y$ is Lipschitz by Theorem~\ref{dm306}, the function $s\mapsto
L(y(s),y'(s))$ is bounded. By the Lebesgue theorem (\ref{eq500})
implies that $\gamma'(s)=-L(y(s),y'(s))$ for a.e.\ $s\in [0,t]$.
The conclusion follows now from (\ref{eq502}) and (\ref{eq503}).
  $\; \; \Box$

\begin{Theorem}\label{dm391b} Let $L\colon \R^n \times \R^n
\rightarrow \R_+$ be a Borel function
and
let $ \varphi \colon \R^n \to \R_+ \cup \{ + \infty \}$
be a lower semicontinuous function with
$\varphi \not\equiv +\infty$. Assume that $({\bf H1})$ and $({\bf
H2})$ are satisfied, and let $V$ be the
value function given by (\ref{value}).
Define the set-valued maps
\begin{eqnarray*}
& F(t,x):=\{u\in\R^n \;| \;  D _{\uparrow}V
(t,x) (-1,u)  \leq -L (x,u)\},
\\
& G(t,x):=\{u\in\R^n \;| \;  D _{\downarrow}V
(t,x) (1,-u)  \geq L (x,u)\}.
\end{eqnarray*}
Given $t>0$, $x\in\R^n$, and
$y\in W^{1,1}(0,t;\R^n)$ with $y(0)=x$, the following conditions are
equivalent:
\begin{itemize}
\item[(a)] $y$ is a
minimizer of (\ref{value}),
\item[(b)] $y' (s) \in
F(t-s,y(s))$ for a.e. $s \in [0,t]$,
\item[(c)] $y' (s) \in
G(t-s,y(s))$ for a.e. $s \in [0,t]$.
\end{itemize}
\end{Theorem}

\noindent
{\bf Proof}  --- \hspace{ 2 mm}
If $y$
is a minimizer if (\ref{value}), then $y'(s) \in F(t-s,y(s))$ and
$y'(s) \in G(t-s,y(s))$ for
a.e.\ $s \in (0,t)$ by (\ref{dm360x}) and (\ref{dm360y}).

If $y$ solves the differential inclusion $y' (s) \in
F(t-s,y(s))$ for a.e. $s \in [0,t]$, we
define
$\gamma (s) :=V(t-s,y(s))$. Since $V$ is locally Lipschitz on
$\R^\star_+\times\R^n$, the function $\gamma$ is locally
absolutely continuous on $[0,t)$. Using (\ref{eq502})
and the definition of $F(t-s,y(s))$, we obtain that
$L(y(s),y'(s)) \leq -\gamma '(s)$ for almost all $s \in [0,t]$,
which implies that $L (y(s),y'(s))$ is integrable on $[0,t-
\varepsilon]$ for every $\varepsilon>0$. By integrating we obtain
\begin{displaymath}
\gamma (t- \varepsilon ) +
\int_{0}^{t- \varepsilon}L (y(s),y'(s))ds \leq \gamma(0).
\end{displaymath}
As $\gamma (t- \varepsilon )=V(\varepsilon,y(t- \varepsilon))$ and
$\gamma(0)=V(t,x)$,
taking the lower limit as $\varepsilon\to 0$ and using Proposition~\ref{prop15}
we get
\begin{displaymath}
\varphi (y(t))+ \int _{0}^{t} L (y(s),y'(s))ds \leq
V(t,x).
\end{displaymath}
Consequently, $y$ is a minimizer of (\ref{value}).

If $y$ solves the differential inclusion $y' (s) \in
G({t-s},y(s))$ for a.e. $s \in [0,t]$,  we repeat the same proof,
replacing (\ref{eq502}) with (\ref{eq503}).
$\;\;\Box$

\vspace{ 5 mm}

\noindent {\bf Acknowledgements.} The work of Gianni Dal Maso is
part of the European Research Training Network ``Homogenization
and Multiple Scales'' under contract HPRN-2000-00109, and of the
Project ``Calculus of Variations'', supported by SISSA and by the
Italian Ministry of Education, University, and Research. Gianni
Dal Maso wishes to thank the hospitality and support of the
University of Paris-Dauphine for two visits in May 2000 and
January 2001, during which the main results of this paper have
been obtained.

The work of H\'el\`ene Frankowska is part of the
European Research Training Network ``Evolution Equations for
Deterministic and Stochastic Control Systems'', RTN-2002-12.
H\'el\`ene Frankowska wishes to thank SISSA for the hospitality in the
summer 2002.


\small
\begin{thebibliography}{99}
\bibitem {amar} AMAR M., BELLETTINI G. \& VENTURINI S. (1998) {\it
Integral representation of functionals defined on curves of\/}
$W^{1,p}$, Proc. Roy. Soc. Edinburgh Sect. A {\bf 128}, 193-217.

\bibitem {ambr} AMBROSIO L.,  ASCENZI O. \&  BUTTAZZO G. (1989) {\it
Lipschitz regularity for minimizers of integral functionals with highly
discontinuous integrands\/}, J. Math. Anal. Appl. {\bf 142}, 301-316.

\bibitem{aubconvex} AUBIN J.-P. (1993)
{\sc Optima and Equilibria},
Grad. Texts in Math. 140, Springer-Verlag, Berlin.

\bibitem{af90sva} AUBIN J.-P. \& FRANKOWSKA H. (1990) {\sc Set-Valued
Analysis}, Birk\-h\"{a}u\-ser, Boston.

\bibitem{bm85} BALL J. \& MIZEL V.J. (1990) {\it One dimensional
variational problems whose minimizers do not satisfy the
Euler-Lagrange equation}, Arch. Rational. Mech. Anal. {\bf 90} (1985),
325-388.

\bibitem{cf91cha}  CANNARSA P. \& FRANKOWSKA H. (1991) {\it Some
characterizations of optimal trajectories in control theory}, SIAM
J. on Control and Optimization, 29, 1322-1347.

\bibitem{cava}   CASTAING C. \& VALADIER M.   (1977)
{\sc Convex Analysis
and Measurable Mutifunctions}, Springer-Verlag, Berlin.

\bibitem{ces83}   CESARI   L.   (1983)   {\sc  Optimization   Theory   and
Applications. Problems with Ordinary Differential Equations},
Appl. Math. 17, Springer-Verlag, Berlin.

\bibitem{Cla} CLARKE F.H. (1983)
{\sc  Optimization   and Nonsmooth Analysis}, Wiley-Inter\-science, New
York.

\bibitem{cla-vin} CLARKE F.H. \& VINTER R.B. (1985)
{\it Regularity properties of solutions to the basic problem in the
calculus of variations\/}, Trans. Amer.
Math. Soc.
{\bf 289}, 73-98.


\bibitem {dm-fra} DAL MASO G. \& FRANKOWSKA (2000) {\it
Value functions for Bolza problems with discontinuous Lagrangians
and Hamilton-Jacobi inequalities\/}, ESAIM Control Optim. Calc.
Var. {\bf 5},  369-394.


\bibitem {DG-But-DM} DE GIORGI E., BUTTAZZO G. \& DAL MASO G. (1983)
{\it On the lower semicontinuity of certain integral functionals\/},
Atti Accad. Naz. Lincei Rend. Cl. Sci. Fis. Mat. Natur. (8) {\bf 74},
274-282.


\bibitem{Ioffe} IOFFE A.D. (1977)
{\it On lower semicontinuity of integral functionals\/},
SIAM J. Control Optim.  {\bf 15}, 521-521 and 991-1000.

\bibitem {Let} LETTA G. (1976)
{\sc  Teoria Elementare dell'Integrazione},
Boringhieri, Torino.

\bibitem{Olech} OLECH C. (1976)
{\it Weak lower semicontinuity of integral functionals\/}, J.
Optim. Theory Appl. {\bf 19}, 3-16.

\bibitem {Roy} ROYDEN H.L. (1969)
{\sc  Real Analysis},
Collier Macmillan, New York.

\bibitem{t1} TONELLI L. (1915) {\it Sur une m\'eth\^ode directe
du calcul des variations}, Rend. Circ. Mat. Palermo {\bf 39}, 223-264.

\bibitem{t2} TONELLI L. (1921) {\sc Fondamenti di Calcolo delle Variazioni},
Vol. 1, 2, Zanichelli, Bologna.

\bibitem{Zie} ZIEMER W.P. (1989) {\sc Weakly Differentiable Functions},
Springer-Verlag, Berlin.
\end{thebibliography}
\normalsize
\end{document}